\documentclass{article}

\usepackage{amsmath, amssymb}

\newtheorem{theorem}{Theorem}[section]

\newtheorem{lemma}{Lemma}[section]
\newtheorem{proposition}{Proposition}[section]
\newtheorem{definition}{Definition}[section]
\newtheorem{remark}{Remark}[section]

\newcommand{\RR}{\mathbb R}
\newcommand{\re}{\hbox{Re }}

\def\dbar{\overline{\partial}}
\def\Re{\hbox{\rm Re}\,}

\def\Aut{\hbox{Aut}\,}

\newfam\msbfam
\font\tenmsb=msbm10   \textfont\msbfam=\tenmsb
\font\sevenmsb=msbm7  \scriptfont\msbfam=\sevenmsb
\font\fivemsb=msbm5   \scriptscriptfont\msbfam=\fivemsb
\def\Bbb{\fam\msbfam \tenmsb}

\def\RR{{\Bbb R}}
\def\CC{{\Bbb C}}
\def\QQ{{\Bbb Q}}
\def\NN{{\Bbb N}}
\def\ZZ{{\Bbb Z}}
\def\II{{\Bbb I}}
\def\TT{{\Bbb T}}
\def\HH{{\Bbb H}}
\def\DD{{\Bbb D}}
\def\BB{{\Bbb B}}
\def\11{{\Bbb 1}}

 \def\HollowBox #1#2{{\dimen0=#1 \advance\dimen0 by -#2
       \dimen1=#1 \advance\dimen1 by #2
        \vrule height #1 depth #2 width #2
        \vrule height 0pt depth #2 width #1
        \llap{\vrule height #1 depth -\dimen0 width \dimen1}%
       \hskip -#2
       \vrule height #1 depth #2 width #2}}
 \def\BoxOpTwo{\mathord{\HollowBox{6pt}{.4pt}}\;}

\def\endpf{\hfill $\BoxOpTwo$}

\begin{document}

\begin{center}
\Large \bf Semi-continuity of Automorphism Groups of Strongly Pseudoconvex Domains: the Low Differentiability Case
\medskip \bigskip \\
\normalsize \rm R. E. Greene, K.-T. Kim, S. G. Krantz, A.-R. Seo
\end{center}

\begin{quote}
{\bf Abstract:}  We study the semicontinuity of automorphism groups
for perturbations of domains in complex space or in complex manifolds.
We provide a new approach to the study of such results for domains
having minimal boundary smoothness.  The emphasis in this study
is on the low differentiability assumption and the new methodology developed accordingly.
\end{quote}

\section{Introduction}

\markboth{R. E. Greene, K.-T. Kim, S. G. Krantz, A.-R. Seo}{Semi-continuity of Automorphism Groups}

It is a familiar perception of everyday life that symmetry is hard to create,
but easily destroyed. To make the crooked straight requires some
definite effort, but the slightest change can suffice to make the
straight a little crooked and hence not straight at all. This perception is easily
substantiated in precise form for geometric objects in Euclidean space.
And it is natural to ask if something similar might apply for automorphism
groups in complex analysis, that is, for the group of biholomorphic self-maps of,
say, a bounded domain in complex Euclidean space.

In one complex variable, this idea does not yield much, at least in the
topologically trivial case. Since all bounded domains that are topologically
equivalent to the unit disc are biholomorphic to the unit disc (Riemann
Mapping Theorem, of course), there is not much interest in discussing how
the automorphism group varies with the domain: it does not vary at all.

But, in higher dimensions, the idea comes into its own. Domains near the unit
ball can have no automorphisms whatever except the identity, and indeed
domains with trivial automorphism group are dense in the set of $C^\infty$
strongly pseudoconvex domains in the $C^\infty$ topology (cf.\
\cite{Greene/Krantz1982} for detailed references to the literature): the
proof of this in fact goes back really to Poincar\'e, in effect, since it depends
essentially only on counting parameters rather than on the details of local
invariant theory, at least once one knows that bihiolomorphic maps extend
smoothly to the boundary \cite{Fefferman}. It is also the case that
domains near the unit ball have automorphism groups which are isomorphic
to a subgroup of the automorphism group of the ball.  Indeed, if a domain is
$C^\infty$ close enough to the ball, the domain is either biholomorphic to the
ball or its automorphism group is isomorphic to a (closed) subgroup of the
unitary group (\cite{Greene/Krantz1982}).

This kind of semicontinuity holds in greater generality
(\cite{Greene/Krantz1982}). If a $C^\infty$ strongly
pseudoconvex domain is not biholomorphic to the ball, then there is a
neighborhood of the domain in the $C^\infty$ topology on the set of all
$C^\infty$ bounded domains with the property that the automorphism
group of every domain in the neighborhood is isomorphic to a subgroup of
the automorphism group of the original domain. (The case of the fixed
domain being biholomorphic to the ball is as in the previous paragraph).

The goal of this paper is to explore the possibility of reducing the level of
differentiability required for this type of result, both for the fixed domain
itself and for the varied domains and the topology upon them. We shall
show in fact that $C^\infty$ can be reduced to $C^2$. This is optimal in the
sense that $C^2$ is the natural setting for the discussion of strong
pseudoconvexity and is the lowest level of regularity for which
the definition is naturally given. (One can of course construct somewhat more
intricate and to some extent artificial ideas of strong pseudoconvexity
wherein the boundary need not have that much regularity, but these will not
be explored here).

It will turn out that the particular complex analysis just discussed can in fact be treated by changing the
whole context to manifolds and general group actions.  The role of complex analysis becomes simply
to guarantee a kind of uniform compactness discussed in Section 2 in detail and in general terms,
momentarily.

To put this matter in perspective, it is desirable to recall in outline how the
semicontinuity results in \cite{Greene/Krantz1982} were obtained.
The starting point is the use of normal family arguments. In this context, the
set-up is as follows. Fix a bounded domain $\Omega_0$. Then a
sequence of bounded domains $\Omega_j$ is considered to converge to
$\Omega_0$ if there is a sequence of maps $\Phi_j\colon \Omega_0 \to
\Omega_j$ which converges to the identity in some
appropriate topology. Now, in this situation, a sequence of automorphisms
$f_j \colon \Omega_j \to \Omega_j$ always has a subsequence $f_{j_k}$
such that the maps $\Phi_{j_k}^{-1} \circ f_{j_k} \circ \Phi_{j_k}$ converge
to some map of $\Omega_0$ to the closure of $\Omega_0$. Here
convergence means uniform convergence on compact subsets of $\Omega_0$.

\newfam\msbfam
\font\tenmsb=msbm10   \textfont\msbfam=\tenmsb
\font\sevenmsb=msbm7  \scriptfont\msbfam=\sevenmsb
\font\fivemsb=msbm5   \scriptscriptfont\msbfam=\fivemsb
\def\Bbb{\fam\msbfam \tenmsb}

\def\RR{{\Bbb R}}
\def\CC{{\Bbb C}}
\def\QQ{{\Bbb Q}}
\def\NN{{\Bbb N}}
\def\ZZ{{\Bbb Z}}
\def\II{{\Bbb I}}
\def\TT{{\Bbb T}}
\def\HH{{\Bbb H}}
\def\DD{{\Bbb D}}
\def\BB{{\Bbb B}}
\def\11{{\Bbb 1}}

The closure may in fact be required. For example, if all the $\Omega_j$ are
the same as $\Omega_0$ with the $\Phi_j$ being the identity, then the
sequence $f_j$ could have a limit that had image in the boundary of
$\Omega_0$, a familiar situation in one variable for the unit disc. For
example, $\Phi_j (z) = \big(z -(1-\frac1j)\big)\big/\big( 1 -(1-\frac1j)
z\big)$ on the disc $\Delta$ in $\CC$ would converge to the constant map $-1$.

However, it is relatively easy to show, and is in fact a classical result that, if
the limit mapping is in fact interior, i.e., if its image lies in $\Omega_0$
itself, then in fact that limit is an automorphism of $\Omega_0$. (A detailed proof is
given in \cite{Krantz}). Thus, in trying to relate the automorphisms of the
$\Omega_j$'s to those of $\Omega_0$, one is interested in situations where it
is guaranteed that the family of maps of the sort described always has
``nondegenerate'' limits, that is , the limits are necessarily the maps into
$\Omega_0$ itself, with no boundary points in the image.

A natural first restriction, arising from looking at the examples where all the
domains are the same, is to those domains $\Omega_0$ which have compact
automorphism group. Then the orbits of the group are necessarily compact
and the limit of any sequence of automorphisms which converges
uniformly on compact sets to some limit will necessarily converge to an
interior limit.

As it happens, every strongly pseudoconvex bounded domains that is not
biholomorphic to the ball has a compact automorphism group. This
was proved by B. Wong \cite{Wong} in the mid 1970s and has been much
generalized since, to the point where the result is not only valid for $C^2$ domains
but is localized completely. If a sequence of automorphisms has
the property that, for some interior point the sequence of the images of the point
converge to a $C^2$ strongly pseudoconvex boundary point of a domain in a
general complex manifold, then the domain is biholomorphic to the ball
(\cite{Efimov}, \cite{Gaussier/Kim/Krantz}). This line of
thought makes it natural to consider the whole normal families situation for
bounded strongly pseudoconvex domains that are not biholomorphic to the
ball, which will indeed be the main topic in this paper. However, certain
aspects of the situation can be treated with no pseudoconvexity invoked
at all.  If one simply assumes the relevant kind of nondegeneracy of normal
families as a hypothesis, then a semicontinuity result already follows. This
matter is treated in Section 2.

It is natural to ask when that hypothesis is satisfied; that is, under what
conditions of a more familiar sort the non-degeneracy condition (stably-interior)
that is required in Section 2 is sure to hold.
As we shall see, it in fact always holds under the hypothesis of $C^2$
strong pseudoconvexity of the boundary of $\Omega_0$ ($\Omega_0$ not
biholomorphic to the ball) and the assumption
that the $\Omega_j$ converge to $\Omega_0$ in the $C^2$ topology. How
this arises requires some explanation.

Already, in \cite{Greene/Krantz1982}, it was observed that non-degeneracy
could be established by considering curvature invariants of the
Bergman metric, at least in the $C^\infty$ case. The argument in
outline was as follows: The Bergman metric of a strongly pseudoconvex
domain is complete in the usual sense of Riemannian geometry
(\cite{Diederich}).  The well-known theorem of Lu Qi-Keng asserts that,
if the Bergman metric had constant holomorphic sectional curvature, then the bounded
domain would be biholomorphic to the ball. Thus, if it is assumed to be not
biholomorphic to the ball, then the holomorphic sectional curvature is not constant.
On the other hand, according to a calculation by Klembeck \cite{Klembeck}
using the Fefferman asymptotic expansion of the Bergman kernel, the
holomorphic sectional curvature approaches a negative constant at the
boundary. (In the usual normalization, the constant is $- 4/(n+1)$, where $n$ is
the complex dimension.)  Let $p$ be a point in the interior where some
holomorphic sectional curvature is not $-4/(n+1)$. Then, since the
holomorphic sectional curvature of the Bergman metric is a biholomorphic
invariant, it follows that there is some positive $\epsilon$ such that the
distance to the boundary of the orbit of $p$ under the automorphism group is greater than or
equal to $\epsilon$. This gives a proof of Bun Wong's theorem on the
compactness of the  automorphism group. But, more significantly from our
viewpoint, it was shown in \cite{Greene/Krantz_Advances} that this
$\epsilon$ can be chosen stably with respect to variation of the domain in the
$C^\infty$ topology. This stability was established by combining interior
stability of the Bergman metric with a (not so easily established) stability of
the Fefferman expansion with respect to variation of the domain.

This program worked, but it was tied specifically to the
$C^\infty$ situation, since the Fefferman expansion requires $C^\infty$
boundary (or at least a large, and rather difficult to determine, number of derivatives).

The semicontinuity of automorphism groups in the $C^2$ case will be
obtained in this paper again by using curvature invariants to bound the
distance of orbits from the boundary stably. But the stability of the
asymptotic constancy of holomorphic curvature of the Bergman
metric will be obtained without using the Fefferman expansion, thus avoiding
the need for a large number of derivatives. Instead, the behavior of the
holomorphic sectional curvature of the Bergman metric will be analyzed
using the ``scaling method,'' as explained in Section 3. The possibility
of using the scaling method depends on noting that the holomorphic sectional
curvature can be expressed in terms of a special basis for the Hilbert
space of square integrable holomorphic functions (cf.\ \cite{Greene/Wu} and
\cite{Epstein} for the special basis concept in generality). This means that
one can detour around the rather awkward formulas from Riemannian
geometry that express the curvature tensor as a whole in terms of the metric
and operate instead with more directly accessible aspects of the fundamental
Bergman construction.

The mechanism by which the general normal family hypotheses introduced
in Section 2 yield semi-continuity results in the sense of isomorphism to
subgroups is the application of the corresponding result in Riemannian
geometry for compact Riemannian manifolds, as established
originally by Ebin \cite{Ebin}: If $g_j$ is a sequence of $C^\infty$
Riemannian metrics on a compact manifold $M$ converging in the $C^\infty$
topology to a $C^\infty$ limit $g_0$ then, for all sufficiently large $j$, the
isometry group of $g_j$ is isomorphic to a subgroup of the isometry group of
$g_0$ via an isomorphism obtained by conjugation by a diffeomorphism of
$M$. This result is actually established in \cite{Ebin} for a finite degree of
differentiability in the sense of Sobolev $H^s$ spaces, but the degree of
differentiability depends on the dimension of $M$, as is typical in $H^s$-space
arguments. Moreover, in the decades since Ebin's paper
\cite{Ebin}, there have been alternative approaches developed
and reductions in the number of derivatives needed. These improved
results will be discussed in Section 5.

To relate this result for compact Riemannian manifolds to the noncompact
case of automorphisms of complex domains, one proceeds as
follows: In the nondegenerate normal family situation already indicated (to
be discussed in detail in Section 2), there can be constructed group-invariant
sub-domains by taking sub-level sets of group-invariant exhaustion functions.
The exhaustion functions are obtained by averaging an arbitrary exhaustion
function with respect to the group action, and the invariant sub-level set can
be taken to have smooth boundary by choosing a sublevel set of a noncritical
value of the invariant exhaustion. Since this invariant sublevel set is
strictly interior, that is has compact closure contained in the domain itself, the
sub-domain will have $C^\infty$ boundary. And the group action on it can
thus be extended to the ``double'' of the sub-domain, regarded as a
compact manifold with boundary. In this doubled situation Ebin's result then
applies directly.

In this setup, the regularity of the boundary of the domain itself plays no role.
As soon as one has the non-degenerate normal family situation, via curvature
invariants or otherwise, then all considerations occur strictly inside the
domain where all mappings involved are holomorphic and hence
$C^\infty$.

However, if one wants to extend to the noncompact case the part of Ebin's
result about diffeomorphism conjugation, then the regularity of the boundary
and of the automorphisms up to the boundary becomes involved. In
the $C^\infty$ case, it was shown in \cite{Greene/Krantz1982}
that in fact one could form the double of the domain itself and extend
the group actions to the double, rather than forming the double of an
invariant sub-domain. Thus the analogue of Ebin's diffeomorphism conjugacy
statement was obtained.

In the last section of this paper, a corresponding result
involving diffeomorphism conjugacy will be obtained for
strongly pseudoconvex domains with low boundary regularity.
For technical reasons, the regularity cannot be quite reduced
to the $C^2$ level which will be all that is needed for the
subgroup semi-continuity. It may be possible that
diffeomorphism conjugacy also applies in the $C^2$ case, but
this result cannot be proved by the methods used here.

It is worth noting that the reference \cite{Greene/Krantz1985}
established a version of the semicontinuity theorem for
automorphism groups in the context of $C^2$ convergence. That
paper was an important first step in the program we are
developing here. The role of holomorphic curvature of the Bergman
metric was replaced by the quotient invariant, that is the
Carath\'eodory volume divided by the Kobayashi-Eisenmann volume.
But the curvature methods here are of independent interest, and
the needed stable uniformity of extension of automorphisms is checked
here in more detail.

The present paper is in some respects a natural continuation of
\cite{Greene/Krantz1985}. In \cite{Greene/Krantz1985}, semicontinuity of
automorphism groups in the sense of isomorphism to a subgroup was established
in the $C^2$ strongly pseudoconvex category (with $C^2$ topology). Isomorphism via a conjugating diffeomorphism was established, however, only in the $C^\infty$ category (\cite{Greene/Krantz1982}), with a program outlined briefly in \cite{Greene/Krantz1985} to establish a conjugating diffeomorphism in the $C^k$ case, $k$ finite but (unspecifiedly) large. This latter was to be based on the Ligocka's extension results for biholomorphic maps. In all cases, the property called in this paper ``stably-interior'' was established using not curvature but rather the quotient of the Carath\'eodory and Kobayashi volume forms. This sufficed for the specific purpose, but it was a less geometrically illuminating biholomorphic invariant than is Bergman metric curvature. But at the time, curvature estimates were only able to be derived from the Fefferman expansion and were hence available only in the $C^\infty$ case.

In the intervening quarter of a century(!) various developments made it possible to view the situation both more broadly and more precisely, the latter in the sense of obtaining specific (low) values for the degree of differentiability needed. These developments include more specific estimates of the differentiability needed  in Ebin's theorem (\cite{Guillemin}, \cite{KimY}) and Lempert's extension theorem (\cite{Lem1}), established here in stable form relative to the variation of the domain. Finally, as shown here (cf.\ also \cite{GKK}) the asymptotic constancy of holomorphic sectional curvature can be analyzed by the scaling method, bypassing the Fefferman expansion and hence obviating the need for $C^\infty$, as already noted.  These developments combined make possible a precise completion of the finite-differentiability program begun in \cite{Greene/Krantz1985}, precise in particular in precise $k$ values.

Useful though the Carath\'eodory-Kobayashi volume quotient was in \cite{Greene/Krantz1985}, it is our perception that the more detailed geometric information provided by the curvature analysis here has more potential for future further applications, as well as being, as we see it, geometrically satisfying in its own right. And the stable Lempert extension estimates also seem to us to have potential for further use, as we hope.


\section{Normal Families and General Semicontinuity of Groups of
Mappings}

In this section, some very general results will be discussed about groups of
diffeomorphisms of open sets in Euclidean spaces. The fundamental idea is
that, as far as semi-continuity of the groups is concerned, the noncompact case
can be converted to the compact case. This is, more precisely, true as far as
semi-continuity in the sense of isomorphism to a subgroup is concerned. We
begin with a definition of an appropriate idea of convergence of the open
sets. For convenience, and without any particular loss of generality, we
restrict our attention to connected open sets, i.e., domains.

\begin{definition} \label{Definition_1}	 \rm
A sequence $\Omega_j$ of connected open sets, or domains, in a Euclidean
space $\RR^n$, is said to {\it containment-converge} to a limit domain
$\Omega_0$ if, for every compact subset $K$ of $\Omega_0$,
$K$ is contained in $\Omega_j$ for all sufficiently large  $j$.
\end{definition}

\begin{definition}   \rm
\label{Definition_2}
If the sequence $\{\Omega_j\}$ of domains containment-converges to a
domain $\Omega_0$, then a sequence of $C^\infty$ mappings $f_j\colon
\Omega_j \to \RR^n$ is said to {\it converge $C^\infty$ normally} if, for
each compact subset $K$ of $\Omega_0$, the mappings $f_j$ and their
derivatives of all orders converge uniformly on $K$.
\end{definition}

Note here that the $f_j$ are defined in a neighborhood of $K$,
any compact set $K$, for all $j$
sufficiently large, so that the desired uniform convergence indeed makes
sense.

For our next definition, we recall that there is a metric, to be denoted $g_K$,
on the set of all $C^\infty$ mappings of a neighborhood of a compact subset
$K$ to $\RR^n$ such that convergence in this metric is equivalent to
convergence of the mappings and their derivatives of all orders
uniformly on the compact set $K$. (cf., e.g., \cite{Greene/Krantz_Advances})

\begin{definition} \label{Definition_3}	  \rm
Suppose that $\{\Omega_j\}$ is a sequence of domains which
containment-converges to a domain $\Omega_0$ and also suppose
that, for each $j$, $G_j$ is a group of diffeomorphisms of $\Omega_j$
and that $G_0$ is a group of diffeomorphisms of $\Omega_0$. We say that
the sequence of groups $G_j$ {\it converges normally} to $G_0$ if,
for each compact subset $K$ of $\Omega_0$ and for each $\epsilon >0$,
there is a $j_{\epsilon, K}$ such that, for each $j >j_{\epsilon, K}$ and
each $\phi_j \in G_j$, the mapping $\phi_j\big|_K$ lies within
$g_{{}_K}$-distance $\epsilon$ of some element of $G_0$.
\end{definition}

In case one has not domains, but compact manifolds and compact groups,
then the situation is as follows:

\begin{lemma}[from \cite{Ebin}, cf.\ \cite{KimY} and \cite{GKK}]
\label{Ebin_based}  \sl
If $M$ is a compact manifold and if $G_j$ is a sequence of compact
subgroups of the diffeomorphism group of $M$ {\rm[}in the topology
determined by the metric $\gamma_M${\rm]} such that $G_j$ converges to
the compact subgroup $G_0$ then, for all $j$ sufficiently large, $G_j$ is
isomorphic to a subgroup of $G_0$. Moreover, the isomorphism can be
obtained by conjugation by a diffeomorphism $\phi_j$ and the
$\phi_j$ can be chosen to converge to the identity {\rm[}again in the
topology determined by the metric $\gamma_M${\rm]}.
\end{lemma}

\noindent {\it Proof}: This result is implied by the result
of D.\ Ebin already alluded to
together with a classical result of Lie group theory. Ebin's result in detail is
that, if $\{g_j\}$ is a sequence of Riemannian metrics on a compact manifold
$M$ which converge in the $C^\infty$ sense to a limit metric $g_0$ then, for
all $j$ sufficiently large, the isometry group of $g_j$ is isomorphic to a
subgroup of the isometry group of $g_0$. Now this result is related to the compact
group situation as follows: With the groups $G_j$ and $G_0$ as above, there
is a metric $g_0$ which is invariant under $G_0$. Then,
because the elements of $G_j$, for $j$ large, are close to
elements of $G_0$, the metric $g_0$ is in an obvious sense close to being
invariant under $G_j$. In particular, averaging $g_0$ with respect to the action
of $G_j$ in the usual fashion produces a metric $g_j$ that is close to $g_0$.
The sequence $\{g_j\}$ converges to $g_0$ in the $C^\infty$ sense.
Ebin's theorem then gives that the isometry group of $g_j$, which of
course includes $G_j$, is isomorphic to a subgroup of the isometry
group of $g_0$. But the connection is not quite complete, since the isometry
group of $g_0$ may in fact be larger than the group $G_0$. But this difficulty
can be handled as follows: Part of Ebin's result is that in fact the isomorphism
to a subgroup can be obtained via conjugation by a diffeomorphism which
can be taken to be close to the identity. So, for all $j$ sufficiently large, we
can choose diffeomorphisms $\psi_j$ such that the group $\widehat{G}_j :=
\psi_j \circ G_j \circ \psi_j^{-1}$ is a subgroup of the group
$\hbox{Isom }(g_0)$ of isometries of $g_0$. Now suppose that the
diffeomorphisms $\psi_j$ converge in the metric $\gamma_M$ defined
above to the identity, as they can certainly be chosen to do. Then the groups
$\widehat{G}_j$, which are subgroups of $\hbox{Isom }(g_0)$, converge to
$G_0$ in the Lie group topology of $\hbox{Isom }(g_0)$: for every open set
$U$ containing $G_0$, there is a $j_U$ such that, for any $j>j_U$, every
element of $\widehat{G}_j$ lies in $U$.

Now one can apply this theorem of Montgomery and Samelson
(\cite{Montgomery/Samelson}):
\sl If G is a compact Lie group and H a closed subgroup, then
there is an open neighborhood U of H such that every subgroup of
G lying in U is conjugate to a subgroup of H. \rm

This now gives the result on compact group actions that we were
seeking. \endpf
\medskip \\

Our goal here is to show how to reduce the domain case to the compact
manifold situation described in the Lemma. Specifically, we want to prove
the following proposition:

\begin{proposition}  \sl
Suppose that $\{\Omega_j\}$ is a sequence of bounded domains in $\RR^N$
which containment-converges to $\Omega_0$ in the sense of
Definition \ref{Definition_1} and that, for each $j$, $G_j$ is a compact group
of diffeomorphisms of ${\rm cl}(\Omega_j)$ and that the sequence $\{G_j\}$
converges $C^\infty$ normally to a compact group $G_0$ of
diffeomorphisms of ${\rm cl}(\Omega_0)$ [convergence in the sense of
Definition \ref{Definition_3}].  Here, of course, \ cl \ denotes the closure of
the indicated set.  Then, for all sufficiently large $j$, the
group $G_j$ is isomorphic to a subgroup of $G_0$.
\end{proposition}

The essential tool is to use group-invariant exhaustion functions to find a
smoothly bounded sub-domain of $\Omega_0$ that is taken to itself by each
element of the group $G_0$ and then pass to the ``double'' of these sub-domains
to form a compact manifold. Then one does a similar construction to
nearby $G_j$-invariant sub-domains of $\Omega_j$ and thus attains the
situation of Ebin's Theorem. We now describe this situation in more detail, following
the arguments developed in \cite{Greene/Krantz_Advances}:

\begin{definition} \label{Definition_4}  \rm
A real-valued function $\rho\colon \Omega \to \RR$ on a domain
$\Omega$ is said to be an{\it exhaustion function} if , for every $\alpha \in
\RR$, the set $\rho^{-1}\big((-\infty, \alpha]\big)$ is compact---that is, the
sub-level sets of $\rho$ are compact.
\end{definition}

Exhaustion functions of course always exist on domains and indeed on
manifolds in general. One for (not necessarily bounded) domains that
frequently occurs in complex analysis is
$\max \big(\|z\|^2, -\log \hbox{ dist}\, ($z$, \hbox{the complement of the domain}) \big)$.
Exhaustion functions with special properties play an important role, for instance,
in the study of Stein manifolds; these are of course more difficult to construct.

Now suppose that $G$ is a compact group of diffeomorphisms on a domain
$\Omega$ and suppose that $\rho$ is an exhaustion function on $\Omega$.
Then the function $\widehat{\rho}$ defined by $\widehat\rho (z) :=
\int_G \rho(g(z))~d\lambda(g)$, where $d\lambda$ is the normalized Haar
measure on $G$, is also an exhaustion function, as one easily sees.
This function is $G$-invariant in the sense that $\widehat\rho (g(z))=
\widehat\rho (z)$. Thus its sub-level sets are invariant under the action of
$G$: a given sub-level set is mapped to itself by each element of $G$.

If $\rho$ is $C^\infty$, then $\widehat\rho$ is also $C^\infty$. In this
case, for all sufficiently large $\alpha$, except for a set of measure 0, the
sub-level set $\widehat\rho^{-1} (-\infty, \alpha]$ is a compact $C^\infty$
manifold-with-boundary. This follows from Sard's Theorem: one need only
take $\alpha$ so large that the sub-level set is nonempty and such
that $\alpha$ is a regular value for $\widehat\rho$.

Now we return to the situation of a sequence of compact groups $G_j$
converging in our previous sense to a compact group $G_0$. As in the
general setting above, we choose a $C^\infty$ exhaustion function $\rho_0$
and average it over $G_0$ to get a $G_0$-invariant, $C^\infty$exhaustion
function $\widehat\rho_0$.

Because $G_j$ is defined on $\Omega_j$ while $\widehat\rho_0$ is
defined on $\Omega_0$, we cannot average $\widehat\rho_0$ to make it
$G_j$-invariant.  We can, however, perform the averaging on arbitrary compact subsets.

Specifically, choose $\alpha$ as above, so that ${\widehat\rho}_0^{-1} (-\infty,
\alpha]$ is nonempty and of course is a compact subset of $\Omega_0$. Let
$L$ be a compact subset of $\Omega_0$ which contains
${\widehat\rho_0}^{-1} (-\infty, \alpha]$ in its interior and let $L_1$
be a compact subset of $\Omega_0$ that contains $L$ in its interior.

Because the sequence $G_j$ converges to $G_0$, it follows easily that, for $j$
sufficiently large, the images under $G_j$ of points of $L$ lie in $L_1$. It
then follows in addition that one can average the function
${\rho}_0$ over the action of $G_j$, as in the process of averaging
to construct $\widehat\rho$. Denote this new function on $L$
by ${\widehat\rho}_j$. Note that, because the elements of $G_j$ are, for $j$
large, close to those of $G_0$, the function ${\widehat\rho}_j$ is
$C^\infty$-close (i.e., $\gamma_L$-close) to
${\widehat\rho}_0$ on $L$. In particular,
the sub-level set $L_1 \cap \widehat\rho_j^{-1} (-\infty, \alpha]$ will be, for $j$
sufficiently large, a smooth manifold-with-boundary which is $C^\infty$
close to ${\widehat\rho}_0^{-1} (-\infty, \alpha]$.

In particular, if we choose a regular value $\alpha$ for ${\widehat\rho}_0$
with the sub-level set $M_0 := {\widehat\rho}_0^{-1}(-\infty, \alpha]$
nonempty then, for all $j$ sufficiently large, the sub-level set $M_j :=
{\widehat\rho}_j^{-1} (-\infty, \alpha]$ will be a nonempty $C^\infty$
manifold-with-boundary. Moreover it will be close to
${\widehat\rho}_0^{-1} (-\infty, \alpha]$ in the $C^\infty$ sense.  Namely,
there will be a sequence of diffeomorphisms $\phi_j\colon M_0 \to M_j$
which converges in the $C^\infty$ sense to the identity on $M_0$.

The next step of the proof is to form the doubles of the invariant sub-domains
with smooth boundary and extend the compact group actions to
them. This will make it possible to apply the lemma above to the
present situation.

For this, suppose that $\Omega$ is a domain, $M$ a compact subset that is a
(nonempty) smooth manifold with boundary and $H$ a compact group of
diffeomorphisms of $\Omega$ that maps $M$ to itself. By the usual
averaging process, similar to the construction of the invariant exhaustion
functions as already discussed, there is a Riemannian metric $g$ on
$\Omega$ for which the elements of $H$ act as isometries,
i.e., $H$ is contained in $\hbox{Isom }(g)$. Now the metric $g$ restricted to
$M$ can be modified so as to remain invariant under $H$ while being a
product metric at and near the boundary of $M$ (see [10] for an early instance
of this construction). This modification is
obtained by first noting that, if $N$ is the inward unit normal (relative to
$g$) along the boundary $\partial M$, then there is an $\epsilon>0$
such that the $g$-exponential map $E\colon \partial M \times [0,\epsilon) \to M$ defined by
$E(p,s) = \exp_p \big(s N(p))$ is a diffeomorphism for $|s| < \epsilon$
and moreover $E(p,s)$, $p \in \partial M$, $0 \leq s < \epsilon$, is a
diffeomorphism of manifolds with boundary onto a neighborhood $V$ of
$\partial M$ in $M$. This is the usual tubular neighborhood construction.
Then one obtains a product metric $h$ on the neighborhood of the boundary
as $h=ds^2 + dp^2$, where $dp^2$ is the metric on $\partial M$ and
we push this metric over via $E$ to the neighborhood $V$ of $\partial M$ in
$M$. This is clearly invariant under $H$. Then one can extend this metric
to all of $M$ in an $H$-invariant way, by taking a function $\phi$ on $V$
that depends on $s$ alone and hence is invariant under the $H$-action. This
function is to be $1$ in a neighborhood of $s=0$, and hence as a function on
$M$, is equal to $1$ in a neighborhood of $\partial M$. And it is to be equal to $0$
when $s>\epsilon/2$. Then $\phi h + (1-\phi) g$ will be a metric on $M$ as
desired: it is smooth on all of $M$, is invariant under $H$, and is a product
metric near $\partial M$.

This metric now extends smoothly to be a metric $\widehat h$ on the double
$\widehat M$ of $M$ in an obvious way. And the group $H$ acts on $\widehat M$ as a subgroup of the isometry group of $\widehat h$. This
subgroup of the isometry group of $\widehat h$ will be denoted by
$\widehat H$.

Our construction can clearly be taken to be stable with respect to the original
$H$-invariant metric $g$ on $M$ in the sense that, if $g_1$ is another
$H$-invariant metric on $M$ which is $C^\infty$ close to $g$, then the
corresponding metric ${\widehat h}_1$ on the double $\widehat M$ of $M$
will be $C^\infty$ close to $\widehat h$.

With these ideas in mind, we return to the convergence situation as before.
Namely, we continue to denote by ${\widehat M}_j$ the doubles of the
$G_j$-invariant sub-level sets, and let ${\widehat G}_j$ denote the
extension of the $G_j$. Now, when $j$ is large, there are diffeomorphisms
$\beta_j\colon {\widehat M}_0 \to {\widehat M}_j$ which have the property that the pullback to ${\widehat M}_0$ of the $G_j$-action on $M_j$
via $\beta_j$ converges in the sense of Lemma \ref{Ebin_based} above.

In particular, $G_j$ is then isomorphic to a subgroup of $G_0$, for all
sufficiently large $j$. Note that, as such, these isomorphisms apply not to
$G_j$ itself but to the restriction of $G_j$ to $M_j$. But, since $M_j$ has
nonempty interior, the restriction of $G_j$ to be an action on the
($G_j$-invariant) set $M_j$ is injective: two isometries of a connected
manifold which are equal on a nonempty open set are equal. (This follows
easily by a standard continuation argument.)  Hence the original $G_j$ are
indeed isomorphic to a subgroup of $G_0$ when $j$ is sufficiently large.
Thus the proposition is established. \endpf

\section{Bergman Metric and Curvature with \boldmath $C^2$ Stability Near the Strongly Pseudoconvex Boundary}

\newfam\msbfam
\font\tenmsb=msbm10   \textfont\msbfam=\tenmsb
\font\sevenmsb=msbm7  \scriptfont\msbfam=\sevenmsb
\font\fivemsb=msbm5   \scriptscriptfont\msbfam=\fivemsb
\def\Bbb{\fam\msbfam \tenmsb}

\def\RR{{\Bbb R}}
\def\CC{{\Bbb C}}
\def\QQ{{\Bbb Q}}
\def\NN{{\Bbb N}}
\def\ZZ{{\Bbb Z}}
\def\II{{\Bbb I}}
\def\TT{{\Bbb T}}
\def\HH{{\Bbb H}}
\def\DD{{\Bbb D}}
\def\BB{{\Bbb B}}
\def\11{{\Bbb 1}}

Let $n>1$ throughout this section. Denote by $\mathcal{D}_n$ the collection
of bounded domains in $\CC^n$ with ${C}^2$ smooth, strongly
pseudoconvex boundary, equipped with the $C^2$ topology via the ${C}^2$ topology on
defining functions. The goal of this section is to establish the following result,
which is Klembeck's theorem [16] for domains in $\mathcal{D}_n$, with
$C^2$ stability. In the statement below the notation $S_\Omega (p;\xi)$ denotes
the holomorphic sectional curvature of the Bergman metric of the domain $\Omega$ at
$p$ along the holomorphic section generated by the tangent vector $\xi$.

\begin{theorem} \label{C2-Klembeck}  \sl
Let $\Omega_0 \in \mathcal{D}_n$.  Then, for every $\epsilon>0$,
there exist $\delta > 0$ and an open neighborhood $\mathcal{U}$ of
$\Omega_0$ in $\mathcal{D}_n$ such that, whenever $\Omega
\in \mathcal{U}$,
$$
\sup \biggl \{ \Big|S_\Omega (p;\xi) - \Big(-\frac4{n+1}\Big) \Big|\colon~
\Omega \in \mathcal{U}, \xi \in \CC^n \setminus \{0\} \biggr \} < \epsilon
$$
for any $p \in \Omega$ satisfying
$\hbox{\rm dis}\, (p,\CC^n \setminus \Omega) < \delta$.
\end{theorem}
\noindent {\it Proof.}
It suffices to show that the following {\it cannot} hold:

\begin{itemize}
\it
\item[(\dag)]  $\exists \epsilon_0
> 0$, $\exists \{\Omega_\nu\} \subset \mathcal{D}_n$ such that
$\Omega_\nu \to \widehat\Omega$ in the ${C}^2$ topology as $\nu \to
\infty$ and $\exists$ a sequence $\{p_\nu \in \Omega_\nu\}$ with
$\displaystyle{\lim_{\nu\to\infty} \hbox{dis }(p_\nu,
\partial\Omega_\nu) = 0}$ such that
$$
\Big| S_{\Omega_\nu} (p_\nu, \xi_\nu) + \frac4{n+1} \Big| \ge
\epsilon_0,
$$
for every $\nu$.
\end{itemize}

Let $\widehat\Omega, \Omega_\nu, p_\nu$ be as in Section 1.
Since the goal is to show that
$$
\lim_{\nu\to\infty} \Big| S_{\Omega_\nu} (p_\nu, \xi_\nu) +
\frac4{n+1} \Big| = 0,
$$
we may assume without loss of generality that
$\displaystyle{\lim_{\nu\to\infty} p_\nu}$ exists. Denote this
limit by $\widehat p$. Notice that $\widehat p \in
\partial\widehat\Omega$.
\medskip

Let $q_\nu \in \partial\Omega_\nu$ be the closest boundary point
of $\Omega_\nu$ to $p_\nu$ for every $\nu=1,2,\ldots$. Then
consider a sequence $R_\nu: \CC^n \to \CC^n$ of complex rigid
motions\index{rigid complex motions} (i.e., unitary maps followed
by translations) in $\CC^n$ and another rigid motion $\widehat R$
satisfying:
\begin{itemize}
\item[(1)] $\widehat R (\widehat p) = 0$ and $R_\nu(q_\nu)= 0$ for
every $\nu$.
\item[(2)] $R_\nu(\partial\Omega_\nu)$ for every $\nu$, and
$\widehat R (\partial\widehat\Omega)$ are tangent at $0$ to the
hyperplane defined by $\re z_1 = 0$.
\item[(3)] $\displaystyle{\lim_{\nu\to\infty} \|R_\nu - \widehat R
\|_{{C}^2} = 0}$, where the norm here is the ${C}^2$-norm of
mappings on an open neighborhood of the closure of
$\widehat\Omega$ in $\CC^n$.
\end{itemize}

Notice that $R_\nu(\Omega_\nu)$ converges to $\widehat R
(\widehat\Omega)$ in the ${C}^2$ topology on bounded domains with
smooth boundaries. Therefore, without loss of generality, we may also assume
the following:
\begin{itemize}
\item[($1'$)] $0 \in \partial\widehat\Omega \cap
\bigg(\displaystyle{\bigcap_{\nu=1}^\infty \partial\Omega_\nu}\bigg)$.
\item[($2'$)] $\partial\widehat\Omega$ and $\partial\Omega_\nu$
(for every $\nu=1,2,\ldots$) share the same outward normal vector
${\bf n} = (-1,0,\ldots,0)$ at the origin.
\item[($3'$)] $p_\nu = (r_\nu, 0, \ldots, 0)$ with $r_\nu > 0$ for every $\nu$.
\end{itemize}
\medskip

Now we need the following three lemmas for the proof.  The first is

\begin{lemma}[\cite{Kim-Yu}, cf.\ \hbox{[8, Ch.\ 10]}]  \sl
\label{localization}
There exists an open neighborhood $U$ of the origin in $\CC^n$ such that
$$
\lim_{\nu \to \infty} \sup_{0 \neq\xi\in \CC^n } \ \left| \frac{2- S_{\Omega_\nu \cap U} (p_\nu; \xi)}{2-S_{\Omega_\nu} (p_\nu;
\xi)} - 1 \right| = 0.
$$
\end{lemma}

Notice that this lemma implies: {\it if
$\displaystyle{\lim_{\nu \to \infty} S_{\Omega_\nu \cap U} (p_\nu;
\xi)}$ exists, it will coincide with $\displaystyle{\lim_{\nu \to
\infty} S_{\Omega_\nu} (p_\nu; \xi)}$.}
\medskip  \\

The next two lemmas convert the problem of understanding the
boundary asymptotic behavior of the Bergman curvature to that of the
stability of the Bergman kernel function in the interior under
perturbation of the boundary:

\begin{lemma}[\cite{Kim-Yu}; cf.\ \hbox{\bf [8, Ch.\ 10]}]  \sl
\label{scaling}
Let the sequence $\{(p_\nu; \xi_\nu) \in \Omega_\nu \times
(\CC^n\setminus\{0\}) \}$ be chosen as above.  Let $B^n$ denote
the open unit ball in $\CC^n$. Then there exists a sequence of
injective holomorphic mappings
$\sigma_\nu : \Omega_\nu \cap U \to \CC^n$
satisfying the following properties:
\begin{itemize}
\item[(\romannumeral1)] $\sigma_\nu (p_\nu) = 0$ {\rm(}the origin of
$\CC^n${\rm)}.
\item[(\romannumeral2)] For every $r$ with $0<r<1$, there exists
$N>0$ such that
$$
(1-r) B^n \subset \sigma_\nu (\Omega_\nu \cap U) \subset (1+r) B^n
$$
for every $\nu > N$.
\end{itemize}
\end{lemma}
\noindent {\it Proof of Lemma \ref{scaling}}:
In our case the situation is simple, because all the points in the sequence
$\{p_\nu\}$ under consideration are located on the $\re z_1$-axis.

Let $\rho$ be a constant with $0<\rho<1$, to be chosen later (depending on $r$).  Let
$$
\mathcal{E}_{\rho} = \{ z=(z_1,\ldots,z_n) \in \CC^n \colon
\re z_1 > (1-\rho) (|z_1|^2 + \ldots + |z_n|^2) \}
$$
and
$$
\mathcal{S}_{\rho} = \{ z \in \CC^n \mid \re z_1 > (1+\rho)
(|z_1|^2 + \ldots + |z_n|^2) \}.
$$
\bigskip  \\

Since $\widehat\Omega$ is a domain with $\mathcal{C}^2$ smooth,
strongly pseudoconvex boundary, there exists an open neighborhood
$U$ of the origin in $\CC^n$ and a biholomorphism-into
$\Psi:U \to \CC^n$ such that
$$
\Psi(\widehat\Omega \cap U) = \{z\in \Psi(U) \mid \re z_1 >
|z_1|^2 + \ldots + |z_n|^2 + R_2(z) \},
$$
where $R_2 (z) = o(|z_1|^2 + \ldots + |z_n|^2)$. Let $V=\Psi(U)$.
Shrinking the neighborhood $U$ if necessary, one obtains that
$$
\mathcal{S}_\rho \cap V \subset \Psi(\widehat\Omega \cap U)
\subset \mathcal{E}_\rho.
$$
Because of the $\mathcal{C}^2$ convergence, and by $(1')$--$(3')$,
one
deduces that there exists $N>0$ such that
$$
\mathcal{S}_\rho \cap V \subset \Psi(\Omega_\nu \cap U) \subset
\mathcal{E}_\rho
$$
for every $\nu > N$. Now let $\lambda_\nu \equiv |\Psi(p_\nu)|$
for every $\nu$. Consider the dilatation maps
$$
\Lambda_\nu (z_1,\ldots,z_n) \equiv
\Big(\frac{z_1}{\lambda_\nu},
\frac{z_2}{\sqrt{\lambda_\nu}}, \ldots,
\frac{z_n}{\sqrt{\lambda_\nu}} \Big).
$$
Notice here that the point sequence $\Psi(p_\nu)$ approaches the
origin non-tangentially to the hypersurface defined by $\Re z_1 = 0$,
which is tangent to $\Psi(\partial\widehat\Omega)$ at the origin.

Finally let
$$
\Phi(z_1,\ldots,z_n) = \Big( \frac{z_1-1}{z_1+1},
\frac{2z_2}{z_1+1}, \ldots, \frac{2z_n}{z_1+1} \Big)
$$
and
$$
\sigma_\nu = \Phi \circ \Lambda_\nu \circ \Psi
$$
for every $\nu$.  Notice that the composition for each $\nu$ by the M\"obius
transformation $\Phi$ adjusts $\sigma_\nu (p_\nu)$ to the origin while
preserving the unit ball.  So, there exists an $\eta \in (0,1)$ such that
$\{\sigma_\nu\}$ yields a sequence of holomorphic  maps satisfying the
desired conclusion. \endpf
\medskip  \\

The third and last lemma toward the proof of Theorem
\ref{C2-Klembeck} is as follows:

\begin{lemma}[\cite{Ramadanov}, \cite{Kim-Yu}; cf.\ \hbox{\bf [8, Ch.\ 10]}]
\label{ramadanov}  \sl
Let $D$ be a bounded domain in $\CC^n$ containing the origin $0$.
Let $\{D_\nu\}$ denote a sequence of bounded domains in $\CC^n$ that
satisfies the following convergence condition:
\begin{quote}
\noindent given $\epsilon>0$, there exists $N>0$ such that
$$
(1-\epsilon) D \subset D_\nu \subset (1 + \epsilon) D
$$
for every $\nu > N$.
\end{quote}
Then, for every compact subset $F$ of $D$, the sequence of Bergman
kernel functions $K_{D_\nu}$ of $D_\nu$ converges uniformly to the
Bergman kernel function $K_D$ of $D$ on $F\times F$.
\end{lemma}

Now we return to the proof of Theorem \ref{C2-Klembeck}.
\medskip  \\

Let $q_\nu, \xi_\nu, \widehat\Omega, \Omega_\nu$ be as above. Let
$U$ be an open neighborhood of the origin as in Lemma \ref{localization}.
Taking a subsequence, we may assume that $q_\nu \in \Omega_\nu \cap U$
for every $\nu$. Select $\sigma_\nu$ as in Lemma \ref{scaling}.

Apply Lemma \ref{ramadanov} (a theorem of Ramadanov [23]) to
our setting, with $D_\nu = \sigma_\nu (\Omega_\nu \cap U)$ and
$D = B^n$. The conclusion of Lemma \ref{ramadanov} states that the
sequence $K_{D_\nu}(z,\zeta)$ converges uniformly to $K_D (z,\zeta)$
on $F \times F$. This of course implies that the sequence $K_{D_\nu} (z,\bar
\zeta)$ converges to $K_D (z,\bar\zeta)$. Notice that the functions now
 involved are holomorphic functions in the $z$ and $\zeta$ variables together.
Therefore Cauchy estimates imply that $K_{D_\nu} (z,\zeta)$ converges
uniformly to $K_D (z,\zeta)$ on $F \times F$ in the $C^k$ sense for any positive
integer $k$. Since the holomorphic sectional curvature of the Bergman metric
involves derivatives of the Bergman kernel function up to fourth order,
we may conclude that  $S_{\sigma_\nu (\Omega_\nu \cap U)} (0; \cdot)$
converges uniformly to $S_{B^n} (0; \cdot)$ on $\{\xi \in \CC^n
\colon \|\xi\|=1 \}$.  Notice that the latter is the constant function with
value $-4/(n+1)$.

Combining this result with the localization lemma (Lemma \ref{localization}),
the conversion lemma (Lemma \ref{scaling}) and the fact that every
biholomorphism is an isometry for the Bergman metric, we see that:
\begin{eqnarray*}
-\frac4{n+1} & = & \lim_{\nu \to \infty} S_{\sigma_\nu (\Omega_\nu
\cap U)} (0; d\sigma_\nu\big|_{q_\nu}(\xi_\nu)) \\
& = & \lim_{\nu \to \infty} S_{\sigma_\nu (\Omega_\nu \cap
U)} (\sigma_\nu (q_\nu); d\sigma_\nu\big|_{q_\nu}(\xi_\nu)) \\
& = & \lim_{\nu \to \infty} S_{\Omega_\nu \cap U} (q_\nu; \xi_\nu) \\
& = & \lim_{\nu \to \infty} S_{\Omega_\nu} (q_\nu; \xi_\nu).
\end{eqnarray*}
This completes the proof of Theorem \ref{C2-Klembeck}.
\endpf
\bigskip  \\

\begin{remark}[Completeness of the Bergman metric]  \rm
The Bergman metric of a bounded strongly pseudoconvex domain is known to
be complete (\cite{Diederich}; for the more general case cf.\ \cite{Ohsawa}).  Since
the scaled limit shown in the proof of Lemma \ref{scaling} is the unit ball, a
variation of that proof-argument also yields the same conclusion as
\cite{Diederich} regarding completeness also (see [8, Section 10.1.7]).
\end{remark}


\section{Stable \boldmath $C^k$-Extension of Automorphisms}

The purpose of this section is to establish the stability of the extension theorem
for the automorphisms of a bounded strongly pseudoconvex domain under $C^k$ perturbation for finite $k$.

\subsection{Convergence of Lempert's Representative Map}

Let $X, Y$ be complex Banach spaces. Let $\phi \colon U \to Y$ be a map from
an open subset $U$ of $X$ into $Y$. The map $\phi$ is said to be {\it differentiable} at
$x\in X$, if there exists a bounded linear map $D_x\phi \colon X \to Y$ such
that
$$
\|\phi(x+h)-\phi(x)-\big(D_x\phi\big)(h) \|_Y= o(\|h\|_X)
$$
as $\|h\|_X\to 0$. Let $L(X,Y)$ denote the set of bounded linear maps from
$X$ into $Y$. It is naturally equipped with the operator norm and hence
becomes a Banach space. Then $\phi$ is said to be $C^1$ on $U$ if $D_x\phi$
exists for all $x \in U$ and $D\phi\colon x \in U \mapsto D_x\phi \in L(X,Y)$
is continuous.

It is also well established what it means for $\phi$ to belong to the class $C^k$
(cf., e.g., \cite{Mujica}).  To understand this point, consider the space
 $L(X \times \cdots \times X, Y)$ of bounded $k$-linear maps with values in
$Y$. For an $S \in L(X \times \cdots \times X, Y)$, define its norm as follows:
$$
\| S \|_k = \sup \{ \|S(h_1, \ldots, h_k)\|_Y \colon \|h_1\|_X \le
1, \ldots, \|h_k\|_X \le 1 \}.
$$
One more piece of notation is necessary:  for a $k$-linear map $S$, a
$(k-1)$-linear map $[S](h)$ is defined by
$$
[S](h)(h_1, \ldots, h_{k-1}) := S (h, h_1, \ldots, h_{k-1}).
$$
Now the idea of a map belonging to the class $C^k$ can be defined
inductively:  the map $\phi$ is said to be $C^k$ at $x\in X$, for $k = 1,2,\ldots$, if there
exits a bounded $k$-linear map
$D^k_x\phi \colon \underbrace{X\times \cdots \times X}_k \to Y $
such that
$$
\|D^{k-1}_{x+h}\phi-D^{k-1}_x\phi-[D^k_x\phi](h)\|_{k-1}=o(\|h\|_X)
$$
as $h\to 0$ and $D^k\phi \colon x \in U \mapsto D_x^k\phi \in
L( \underbrace{X\times \cdots \times X}_k, Y)$ is continuous. It
is also known that such a $D^k_x\phi$ is symmetric $k$-linear.

Similarly, we may define the concept of H\"older class. For an $\alpha$
with $0<\alpha \leq 1$, a map $\phi$ is said to belong to the class $C^{k,\alpha}$
if $\phi$ is $C^k$ and
$$
\sup_{x,y\in U \atop x \neq y} \frac{\| D^k_x \phi -
D^k_y \phi\|_{k}}{\|x-y\|_X^\alpha} < \infty.
$$

Throughout this section, we denote by $\Delta$ the open unit disc $\{
z\in \CC \colon |z|<1 \}$.  We shall follow the terminology of
\cite{Lem1} closely.  Let $s$ be such that $0<s<\alpha$ and set
\begin{eqnarray*}
X_n &=& \{ f : \partial\Delta \to \CC^n \mid f \in C^{0,s} \}\\
Y_n &=& \{ f\in X_n \colon f \ \hbox{admits a holomorphic
continuation to} \ {\rm cl}(\Delta) \} \\
Y_n^\perp &=& \{ f\in X_n \colon f
\hbox{ admits an anti-holomorphic continuation} \\
& & \qquad \qquad \qquad\hbox{ to ${\rm cl}(\Delta)$ with } f(0)=0 \}.
\end{eqnarray*}
Notice that $X_n=Y_n\oplus Y_n^\perp$.
\smallskip   \\

Let $\Omega=\Omega_\rho$ be a bounded strictly convex domain defined
by the  $C^{k+1,\alpha}$ defining function $\rho$. Then there exists a convex
open neighborhood $V$ of ${\rm cl}(\Omega)$ such that  $\Omega =
\Omega_\rho = \{ z \in V \colon \rho(z)<0\}$, where the defining function
$\rho\colon  U \to \RR$, defined on a convex open set $U$ with ${\rm cl}(V)
\subset  U$, is of class $C^{k+1, \alpha}$ $(k \geq 1, 0<\alpha<1)$ with $d\rho
\neq 0$ at any point of $\partial\Omega$. We may further assume
without loss of generality that
\begin{itemize}
\item[(1)] $\rho\colon U \to \RR$ is compactly supported

and

\item[(2)] the real Hessian of $\rho$ is strictly positive at every point
of $\partial\Omega$.
\end{itemize}
\vspace*{.12in}

Let $\mathcal{N}$ be a $C^{k+1,\alpha}$ neighborhood of $\rho$
chosen so small that every element of $\mathcal{N}$ has  its real
Hessian strictly positive at every point of $V$.
We may require further that there exists a constant $R'>0$ such that, if
$\eta,\tau\in \mathcal{N}$, then
$\|\eta -\tau\|_{C^{k+1,\alpha}(U)}<1$ and
$\|\eta\|_{C^{k+1,\alpha} (U)}<R'$.

Let $p$ be a point in $\Omega$ and let $W$ a neighborhood of $p$
 in $\Omega$ such that $W\subset \Omega_\eta$ for all
$\eta\in \mathcal{N}$. Define
$
\Theta \colon \mathcal{N}\oplus (\CC^n\setminus \{0\}) \oplus W
\to Y_n
$
by
$\Theta(\eta,\zeta, q)=e_{\eta,\zeta,q} $, where $e_{\eta,\zeta,q}$ is the stationary map
(= extremal map) from ${\rm cl}(\Delta)$ to ${\rm cl}(\Omega_\eta)$ satisfying
$e_{\eta,\zeta,q}(0)=q$ and
${e_{\eta,\zeta,q}}'(0)=\mu  \zeta$ for some $\mu >0$.

\begin{proposition}\label{A}  \sl
The map $\Theta$ is locally $C^{k,\alpha -s}$ for any $0<s<\alpha$.
\end{proposition}
\noindent {\it Proof:}
Let $(\eta,v,q)\in \mathcal{N}\oplus
\CC^n\setminus \lbrace 0\rbrace \oplus \mathcal{W}$.
We shall prove that $\Theta$ is $C^{k,\alpha-s}$ near $(\eta,v,q)$.
Let $e=e_{\eta,v,q}=(e_1,...,e_n):{\rm cl}(\Delta) \to {\rm cl}(\Omega_\eta)$
and $\tilde{e}=(\tilde{e}_1,...,\tilde{e}_n)$ be the {\it dual map} of
$e$. (See \cite{Lem0} for the definition of the dual map and its
basic properties.)  Since $\tilde{e}$ has no zeros, there exist two components
which do not vanish simultaneously by a generic linear change of coordinates.
Hence we may assume without loss of generality that $\tilde{e}_1$ and
$\tilde{e}_2$ do not vanish simultaneously on ${\rm cl}(\Delta)$,.
It is also shown in \cite{Lem0} that $\tilde{e}$ extends to a  $C^{k,\alpha}$
map up to the boundary, and that there exist functions
$G_1, G_2\in C^{k,\alpha}({\rm cl}(\Delta))$ that are holomorphic in $\Delta$ and satisfy
$\tilde{e}_1G_1+\tilde{e}_2 G_2 \equiv 1$.
Define the holomorphic matrix $H$ on $\Delta$ by
$$
H=
\left ( \begin{array}{ccccc}
e_1' & -\tilde{e}_2 & -G_1 \tilde{e}_3 &\cdots& -G_1 \tilde{e}_n\\
e_2' &\tilde{e}_2   & -G_2 \tilde{e}_3 & \cdots& -G_2 \tilde{e}_n\\
e_3' &       0      &         1        & \cdots &          0      \\
\vdots  &   \vdots       &           \vdots    & \ddots &      \vdots     \\
e_n' &       0      &         0        &\cdots&          1
\end{array}
\right )
$$
Notice that $H\in C^{0,s}({\rm cl}(\Delta))$ and $\det(H)\neq 0$ on ${\rm cl}(\Delta)$.
Set
$$
Y_n^{R,U}= \biggl \{ f\in Y_n \colon
\|f\|_{C^1({\rm cl}(\Delta))}<R, f(\partial \Delta)\subset U
\biggr \}
$$
and define the map
$$
\Phi:\mathcal{N}\oplus \CC^n\setminus \lbrace 0 \rbrace
\oplus \mathcal{W} \oplus Y_n^{R,{U}}\oplus \RR
\to
T \oplus Y_{n-1}^\perp \oplus\CC^n\oplus\CC^n
$$
by
$$
\Phi(r,v,q,f,\lambda)= \left ( r \circ f, \pi\Bigl(\frac{\langle H^tr_z\circ
f \rangle}{(H^tr_z\circ f)_1}\Bigr), f(0)-q, f'(0) - \lambda v \right ) \, ,
$$
where:
\begin{itemize}
\item[(\romannumeral 1)]
$T=\{ g :\partial \Delta \to \RR \colon g \in C^{0,s} \}$,
\item[(\romannumeral 2)]
$\pi : Y_{n-1}\to Y_{n-1}^\perp$
is defined by
$\displaystyle{
\pi \bigl(\displaystyle\sum_{-\infty}^{\infty}a_kz^k
\bigr)=\displaystyle\sum_{-\infty}^{-1}a_kz^k
}$, and
\item[(\romannumeral 3)]
$(H^tr_z\circ f)_j$ denotes the $j$-th component of $H^t r_z\circ f$ and \hfill\break
$\langle H^t r_z\circ f\rangle =((H^tr_z\circ f)_2,...,(H^tr_z\circ f)_n)$.
\end{itemize}
Then $f:{\rm cl} (\Delta)\to {\rm cl}(\Omega_r)$ is an extremal map
satisfying $f(0)=q, f'(0)=\lambda v$
if and only if $\Phi(r,v,q,f,\lambda)=0$.
So, according to \cite{Lem1}, we only need to prove that $\Phi$
is $C^{k,\alpha -s}$. For this purpose define the map
$\Psi : \mathcal{N} \oplus Y_n^{R,U} \to T$
by $\Psi(r,f)=r\circ f$. Then we pose the following
\medskip  \\

\noindent {\bf Claim.}  $\Psi$ is $C^{k,\alpha -s}$.
\medskip  \\

We shall prove this claim by induction on $k$.
We need some notation. For a domain $\Omega$, $k \in \ZZ^+$, and $0 < \alpha \leq 1$,
denote by
$$			
 \|g\|_{C^{k,\alpha}({\rm cl}(\Omega))}
=\sup_{\small x \in {\rm cl}(\Omega)\atop
|\gamma|=0,1,..,k}|D^\gamma g(x)|
+ \sup_{\small x,y\in
{\rm cl}(\Omega)\atop x\neq y, |\gamma|=k}
\frac{|D^\gamma g(x)-D^\gamma g(y)|}{|x-y|^\alpha}.$$
Moreover, $A\lesssim B$ will mean that $A\leq CB$ for some constant
$C$.  In turn,
$A \lessapprox B$ will mean that $A \rightarrow 0$ whenever $B
\rightarrow 0$.

Let  $j \in \{0,\ldots, k\}$. Let $\mathcal{N}_j =\lbrace r \in C^{j+1,\alpha}(U)\colon\|r\|_{C^{j+1,\alpha}(U)}<R'\rbrace$.
Define $\Psi_j \colon\mathcal{N}_j \oplus Y_n^{R,U}\to T$
by $\Psi_j(r,f)=r\circ f$.
Suppose that, for all $r,\tau \in \mathcal{N}_j$, we have
$\|r-\tau\|_{C^{j,\alpha}(U)}<1$.

In case $j=0$, it suffices to show that
$$
\|\Psi_0(r,f)-\Psi_0(\tau,g)\|_{C^{0,s}(\partial \Delta)}
\lesssim
\biggl (\|r-\tau\|_{C^{0,\alpha}(U)}+\|f-g\|_{C^{0,s}(\partial
\Delta)} \biggr )^{\alpha -s}
$$
For $x\in \partial \Delta$,
\begin{eqnarray*}
| r\circ f(x)-\tau \circ g(x)|
& \leq & |r\circ f(x)-r\circ g(x)|+|r\circ g(x)-\tau\circ g(x)|\\
& \lesssim & |f(x)-g(x)|^{\alpha -s} + |(r-\tau)\circ g(x)| \\
& \lesssim & \biggl (\|f-g\|_{C^{0,s}(\partial \Delta)}+\|r-\tau
\|_{C^{0,\alpha}(U)}\biggr )^{\alpha -s} \, .
\end{eqnarray*}

For $x,y\in \partial \Delta$, let $\delta(x,y) = r\circ
f(x)-\tau\circ g(x)-r\circ f(y)+\tau\circ g(y)$. Then
\begin{eqnarray*}
|\delta (x,y)|
&\leq& |r\circ f(x)-r\circ g(x)|+|r\circ g(x)-\tau\circ
g(x)|\\
& & \qquad +|r\circ f(y)-r\circ g(y)|+|r\circ g(y)-\tau\circ g(y)|\\
&\leq&  2(R')^\alpha |f(x)-g(x)|^\alpha +
2\|r-\tau\|_{C^{0,\alpha}(U)} \\
&\leq&  2(RR')^\alpha \|f-g\|_{C^{0,s}(\partial \Delta)}^\alpha
+2\|r-\tau\|_{C^{0,\alpha}(U)}
\end{eqnarray*}
and
\begin{eqnarray*}
|\delta(x,y)| &\leq& |r\circ f(x)-r\circ f(y)|+|\tau\circ
g(x)-\tau\circ g(y)|\\
&\leq & R'|f(x)-f(y)|^\alpha +R'|g(x)-g(y)|^\alpha \\
&\leq&  2RR'|x-y|^\alpha.
\end{eqnarray*}
This implies that
$$
|\delta(x,y)|
\lesssim \bigl(\|f-g\|_{C^{0,s}(\partial \Delta)}
+\|r-\tau\|_{C^{0,\alpha}(U)}\bigr)^{\alpha -s}|x-y|^s,
$$
which proves the case $j=0$. \\

Let $j>0$. Suppose that $\Psi_j:\mathcal{N}_j\oplus Y_n^{R,U}\to
T$ is of class $C^{j,\alpha -s}(U)$.
Then, since
$$
D_{(r,f)}\Psi_{j+1}(\tau,g)=(r'\circ f)g+\tau\circ
f=\Psi_j(r',f)g+\Psi_j(\tau,f),
$$
it follows that $\Psi_{j+1}$ is of $C^{j+1,\alpha -s}(U)$. This proves the
claim.\\

Since $\pi$ is a bounded linear map, the second component of
$\Phi$ is also of class $C^{k,\alpha -s}(U)$. The proof of the proposition is now complete.
\endpf
\bigskip  \\

Next, for $r\in \mathcal{N}, q\in W$, consider {\it Lempert's
representation map} at $q$ for the domain $\Omega_r$.  We have $L_{r,q} : {\rm cl}(\BB^n) \to {\rm cl}(\Omega_r)$
defined by
$L_{r,q}(\zeta)=\Theta(r,\zeta, q)(|\zeta|)=e_{r,\zeta,q}(|\zeta|)$. The following proposition discusses the convergence of these representation maps.

\begin{proposition} \label{L}  \sl
Let $\rho_j,\rho \in \mathcal{N}$ and let $p_j,p \in W$ be
such that $\|\rho_j-\rho\|_{C^{k+1,\alpha}(U)} \to
0,|p_j-p| \to 0$ as $j\to \infty$. Set the notation
$L_j:=L_{\rho_j,p_j}$,  $L:=L_{\rho,p}$ and $\BB^n_\delta :=
\BB^n \setminus
\lbrace z\in \CC^n \colon |z|<\delta \rbrace$.
Then, for $0<\beta<\alpha$ and $0<\delta<1$, Lempert's representation maps
$L_j$ for $\Omega_{\rho_j}$ converge to Lempert's representation map $L$ for $\Omega_\rho$ on $\BB^n_\delta$ in the
$C^{k,\beta}$ norm, as $j\to \infty$.
\end{proposition}
\noindent {\it Proof:}
Let $\hbox{ev} \colon Y_n \to \CC^n$ be defined by  $\hbox{ev}(g)=g(1)$
(here ``ev'' stands for ``evaluation'' map).  Since $L(\zeta)=\Theta(\rho,\zeta, p)(1)
=\hbox{ev}\circ\Theta(\rho,\zeta,p)$ for $\zeta\in \partial\BB^n$,
 $\hbox{ev}$ is bounded linear.
Write $D^\ell=\frac{\partial^{m_1+...+m_n}}{\partial
x_1^{m_1}...\partial x_n^{m_n}}$, where $|\ell|=m_1+...+m_n$. Then
$D^\ell L(\zeta)=(D_{(\rho,\zeta,p)}^{|\ell|}\Theta)
(\underbrace{\vec{x_1},...,\vec{x_1}}_{m_1};...;\underbrace{\vec{x_n},...,\vec{x_n}}_{m_n})(1)$.
So $\|L_j-L\|_{C^{k,\beta}(\partial \BB^n)} \to 0$ as
$j\to \infty$.
\smallskip \\

Given  $v\in \CC^n,|v|=1, \xi\in \Delta$, denote by $e$ the extremal map
satisfying $e(0)=p, e'(0)=\mu v$ for some $\mu >0$. Then
$L(\xi v)=e(\xi |v|)=e(\xi)$.  This implies that $L(\xi v)$ is holomorphic
with respect to $\xi$.  Now the Poisson integral formula for $\Delta$ yields the desired conclusion.
\endpf
\smallskip \\

\subsection{A Simultaneous Extension Theorem for Automorphisms}
The next goal is to establish the following theorem:

\begin{theorem}[Uniform extension]  \sl
Let $\Omega_j, \Omega$ be strongly pseudoconvex, \\
bounded domains
in $\CC^n$
with $C^{k+1,\alpha} ( k \in \ZZ , k\geq2 , 0<\alpha \leq 1)$
boundaries such that $\Omega_j$ converges to $\Omega$
as j $\rightarrow \infty$ in the $C^{k+1,\alpha}$ topology.
Let a sequence $\{f_j \in \Aut(\Omega_j) \colon j=1,2,\ldots\}$ be given.
Then, for any $\beta$ with  $0<\beta < \alpha$, the sequence $f_j$ {\rm (}every one of which extends to a $C^{k, \beta}$ diffeomorphism of the closure ${\rm cl}(\Omega_j)$ by the `sharp extension theorem' of Lempert \cite{Lem1}{\rm )}  admits a subsequence $\Omega_{j_\ell}$ and $f_{j_\ell} \in \Aut(\Omega_{f_\ell})$ that converges to the  $C^{k, \beta}$-diffeomorphism, the  extension of $f \in \Aut(\Omega) $, in the $C^{k, \beta}$ topology.
\label{SET}
\end{theorem}

This indeed is a normal family theorem together with
H\"older convergence up to the boundary.  Of course precise definitions and terminology are in order, which will be presented here as the exposition progresses.

\begin{definition}  \rm
Let $\Omega_j$ and $\Omega$ be bounded strongly pseudoconvex domains in
$\CC^n$ with $C^{k,\alpha} ( k \in \ZZ , k\geq2 , 0<\alpha \leq 1)$
boundaries.  As $j \to \infty$, the sequence of domains $\Omega_j$ is said to {\it converges to $\Omega$ in
the $C^{k,\alpha}$ topology}, if there exist an open neighborhood $U$ of
${\rm cl}(\Omega)$, $C^{k,\alpha}$ diffeomorphisms $ F_j: U \rightarrow
U$, and a positive integer $N$ such that:
\begin{itemize}
\item ${\rm cl}(\Omega) \subset\subset U$, %
\item ${\rm cl}(\Omega_j) \subset\subset U $ for all $j>N $,
\item each $F_j $ maps  ${\rm cl}(\Omega)$ onto ${\rm cl}(\Omega_j)$ as a
$C^{k,\alpha}$ diffeomorphisms, for every $j>N$, and
\item $\|F_j-\hbox{\rm id}\|_{C^{k,\alpha}(U)} \rightarrow 0$ and
$ \|F_j^{-1}-\hbox{\rm id}\|_{C^{k,\alpha}(U)} \rightarrow 0 $, as $j \rightarrow \infty$.
\end{itemize}

In a similar manner, we say that the sequence of maps $f_j\in C^{k,\alpha}(\Omega_j,\CC^m)$ {\it converges to $f\in C^{k,\alpha}(\Omega, \CC^m)$ in the $C^{k,\alpha}$ sense}, if  $\lim_{j\to\infty}\|f_j\circ F_j - f\|_{C^{k,\alpha}(\Omega)} = 0$.
\end{definition}

We now present several technical lemmas.

\begin{lemma} \label{T}  \sl
Let $\Omega_j$ be a domain in $\RR^{n_j}$ for each  $j=1,2,3$, respectively.  If
\begin{itemize}
\item[(\romannumeral 1)] $g, h\colon \Omega_1 \to \Omega_2$ are
$C^{k',\alpha'}$ maps that are injective,
\item[(\romannumeral 2)] $f \colon \Omega_2 \to \Omega_3$ is a
$C^{k'', \alpha''}$ map,

and

\item[(\romannumeral 3)]  $(k,\alpha)$ is the pair of the positive integer $k$ and the real number $\alpha$ satisfying $k+\alpha=\min \{k'+\alpha',k''+\alpha''\}$ and  $0<\alpha \leq 1$,
\end{itemize}
then
\begin{itemize}
\item[(1)] $f \circ g \in C^{k,\alpha}(\Omega_1, \Omega_3)$

and

\item[(2)] $\|f\circ g-f\circ h\|_{C^{k,\beta}(\Omega_1)}\lessapprox
\|g-h\|_{C^{k,\alpha}(\Omega_1)}$ for any $\beta$ with  $0< \beta <\alpha$.
\end{itemize}
\end{lemma}
\noindent {\it Proof:}
We shall present the verification of (1) only, as our arguments are mostly by straightforward computation and the proof of (2) is similar.
The chain rule implies that
$$
D^\ell (f\circ g)(x) =
\sum (D^m
f)(g(x))(D^{m_1}g(x))^{m_1'}(D^{m_2}g(x))^{m_2'}...(D^{m_n}g(x))^{m_n'}	  \, ,
$$
where $\ell$ and $m$ are multi-indices and $ m_j$ nonnegative integers satisfying
$|m|\leq |\ell|$ and $\sum m_j'\leq |\ell|$. (We use the usual
multi-index notation here; we omit detailed expressions as they are standard.)
Note that
$$
 \|f\circ g\|_{C^{k,\alpha}}
=\sup_{\small x\in \Omega_1 \atop
\small 0 \le |\gamma|\le k}|D^\gamma (f\circ g)(x)|
+ \sup_{\small x,y\in \Omega_1 \atop x\neq y,
|\gamma|=k}
\frac{|D^\gamma (f\circ g)(x)-D^\gamma (f \circ g)(y)|}{|x-y|^\alpha}.
$$
First, one sees immediately that
$$
\sup_{\small \substack x\in \Omega_1 \atop
|\gamma|=0,1,..,k}|D^\gamma (f\circ g)(x)|
\lesssim \|f\|_{C^{k,\alpha}(\Omega_2)} \sum
\|g\|^{m_1'+\ldots +m_n'}_{C^{k,\alpha}(\Omega_1) }< \infty.
$$
On the other hand,
\begin{eqnarray*}
& & |D^\gamma (f\circ g)(x) - D^\gamma  (f \circ g)(y)|  \\
&=&\bigg|\sum \Big\{ D^m f(g(x))\cdot (D^{m_1}g(x))^{m_1'} \cdot \ldots \cdot (D^{m_n}g(x))^{m_n'}
 \\
& & \qquad\qquad -D^m f(g(y))\cdot (D^{m_1}g(y))^{m_1'} \cdot \ldots \cdot (D^{m_n}g(y))^{m_n'}\Big\}\bigg|
\\
& \le &\sum \bigg\{\Big|\big(D^m f(g(x))-D^m
f(g(y))\big) \cdot \big(D^{m_1}g(x))^{m_1'}\cdot \ldots \cdot (D^{m_n}g(x))^{m_n'} \Big|
\\
&  & \quad+  \Big|\big(D^m f(g(y))\big) \cdot \big( (D^{m_1}g(x))^{m_1'} - D^{m_1}g(y))^{m_1'}\big) \\
& & \qquad\qquad\qquad\qquad \qquad\qquad\cdot (D^{m_2}g(x))^{m_2'}\cdot \ldots \cdot(D^{m_n}g(x))^{m_n'} \Big|
\\
& & \quad +
\\
& & \quad \vdots
\\
& & \quad + \ \Big|\big(D^m f(g(y))\big) \cdot \big(D^{m_1}g(y))^{m_1'}\big)
\cdot \ldots \cdot \big( (D^{m_n}g(x))^{m_n'}-(D^{m_1}g(y))^{m_1'} \big) \Big| \bigg\}
\\
& &\lesssim
\|f\|_{C^{k,\alpha}(\Omega_2)}(1+\|g\|_{C^0(\Omega_1)}^\alpha)P(\|g\|_{C^{k,\alpha}(\Omega_1)})
|x-y|^\alpha,
\end{eqnarray*}
where $P$ is an appropriate polynomial with $P(0,...0)=0$. Hence (1) follows. We omit the proof of (2).
\endpf
\smallskip \\

\begin{lemma}  \sl
\label{tt}
Let $k \geq 1$. Assume that $\Omega_1, \Omega_2$ are bounded domains in $\RR^n$ admitting $C^{k,\alpha}$ diffeomorphisms $f_j, f \colon  {\rm cl}(\Omega_1) \to {\rm cl}(\Omega_2)$
satisfying $\|f_j-f\|_{C^{k,\alpha}({\rm cl}(\Omega_1))}
\rightarrow 0$ as $j \rightarrow \infty $.
If $\displaystyle \lim_{j\to \infty} \sup_{x\in {\rm cl}(\Omega_2)}
|f_j^{-1}(x)-f^{-1}(x)| = 0$,
then $\displaystyle \lim_{j\to \infty} \|f_j^{-1}-f^{-1}\|_{C^{k,\beta}({\rm cl}(\Omega_2))}
= 0$ for any
$0<\beta<\alpha$.
\end{lemma}
\noindent {\it Proof:}
The inverse function theorem implies that
$df^{-1}_j \big|_{f_j(y)}=\big(df_j\big|_y\big)^{-1}$ and \newline
$df^{-1} \big|_{f(y)}=\big(df\big|_y\big)^{-1}$.
Since ${\rm cl}(\Omega_1)$ and ${\rm cl}(\Omega_2)$ are compact,
there exist a constant $C>0$ and a positive integer  N such that
$|\det(df|_y)|>C$ and $|\det(df_j|_y)|>C$ for any point
$y\in\Omega_1$ and any integer $j>N$.  Lemma $\ref{T}$ and its proof-argument above now yield the desired conclusion.
\endpf
\smallskip \\

\begin{lemma}   \sl
\label{TT}
Let $k$ be an integer with $k \ge 2$ and $\alpha$ a real number satisfying
$0< \alpha \leq 1$.  If $\Omega$ is a bounded, strongly pseudoconvex domain
in $\CC^n$ with $C^{k+1,\alpha}$ boundary then, for any $\beta$ with  $0< \beta <\alpha$, there exist an open neighborhood $\mathcal{U}$ of $\Omega$ and a constant $C$ such that $\|f\|_{C^{k,\beta}({\rm cl}(\Omega'))}<C$
for any $\Omega' \in \mathcal{U}$ and any $f\in \Aut(\Omega')$.
\end{lemma}
\noindent {\it Proof:}
Assume the contrary. Then there exists a sequence of strongly pseudoconvex domains $\Omega_j$ with $C^{k+1,\alpha}$ boundary converging to $\Omega$ in the $C^{k+1,\alpha}$ topology and a sequence $f_j \in \Aut(\Omega_j)$ such that
$$
\lim_{j\to \infty }
\|f_j\|_{C^{k,\beta}({\rm cl}(\Omega_j))}=\infty.
$$
Then either
\begin{itemize}
\item[(1)] there exists a sequence $\{x_j \in \Omega_j\colon j=1,2,\ldots\}$ such that $|D^\gamma f_j(x_j)| \rightarrow \infty$
as $j\rightarrow \infty$ for some multi-index $\gamma$ satisfying $0\leq|\gamma|\leq
k$;
\end{itemize}
or
\begin{itemize}
\item[(2)] there exist $x_j,y_j \in \Omega_j$ such that
$\displaystyle{ \lim_{j\rightarrow \infty} \frac{|D^\gamma f_j(x_j) -D^\gamma f_j(y_j)|}{|x_j-y_j|^\beta} =\infty}$
for some multi-index $\gamma, |\gamma|=k$.
\end{itemize}
\vspace*{.15in}

Suppose that  (1) holds.  Then, since the sequence $f_j$ converges to $f$ in the $C^\infty (K)$ topology
on every compact subset $K$ of $\Omega$, it must be the case that $\lim_{j\to\infty} x_j = p\in \partial \Omega$ (taking a subsequence if necessary).

We shall arrive at the desired contradiction to (1) by means of the following three steps:
\bigskip  \\

\noindent
\bf Step 1. \it Adjustments. \rm
\medskip  \\

Here $F_j$ denotes the same diffeomorphism of $\hbox{cl}(\Omega)$ onto $\hbox{cl}(\Omega_j)$ as in Definition 4.1.
Set $F_j(p)=p_j,\; f_j(p_j)=q_j ,\; f(p)=q$. Take the invertible affine
$\CC$-linear transformations $T, T_j, t,t_j\colon\CC^n\to \CC^n$
such that
\begin{itemize}
\item  $T_j(p_j)=T(p)=t_j(q_j)=t(q)=(0,\ldots,0)$;
\item the outward normal vectors to the boundary of $T_j(\Omega_j), T(\Omega), t_j(\Omega_j)$ and $t(\Omega)$, respectively, at $(0,\ldots,0)$
are equal to $(1,0,...,0)$; and
\item  $\displaystyle{\lim_{j\to\infty} T_j = T}$ and $\displaystyle{\lim_{j\to\infty} t_j = t}$.
\end{itemize}

Then  $T_j(\Omega_j)$ converges to $T(\Omega)$ in the
$C^{k+1,\alpha}$ topology, and also $t_j(\Omega_j)$ converges to $t(\Omega)$.
Replacing therefore $f$ and $f_j$, respectively, by $t\circ f\circ T^{-1}$ and
$t_j \circ f_j\circ T_j^{-1}$, we may assume that
\begin{itemize}
\item $\Omega, \Omega_j, \widehat{\Omega}, \widehat{\Omega}_j$
are bounded strongly pseudoconvex domains
with $C^{k+1, \alpha}$ boundaries
such that $\Omega_j$ (and $\widehat{\Omega}_j$, respectively)
converges to $\Omega$ (and to  $\widehat{\Omega}$, respectively) in the
$C^{k+1,\alpha}$ topology.  More precisely, there exist a neighborhood ${U}$ (and  $\widehat{U}$, respectively) of ${\rm cl}(\Omega)$ (and of ${\rm cl}(\widehat{\Omega})$, respectively) and diffeomorphisms
$F_j \colon {\rm cl}(\Omega) \to {\rm cl}(\Omega_j)$ and $\widehat{F}_j \colon {\rm cl}(\widehat\Omega) \to {\rm cl}(\widehat\Omega_j)$
such that
$F_j(0)=\widehat{F}_j(0)=0$ and
the maps $F_j, F_j^{-1}, \widehat{F_j}$ and $\widehat{F_j}^{-1}$converge to
the identity map in the $C^{k+1,\alpha}$ sense.
\item $\rho, \rho_j=\rho\circ F_j^{-1}, \widehat{\rho},
\widehat{\rho}_j=\widehat{\rho}\circ \widehat{F}_j^{-1}$
are defining functions of
$\Omega,\Omega_j,\widehat{\Omega},\widehat{\Omega}_j$,
respectively, such that
$\|\rho-\rho_j\|_{C^{k+1,\beta}({U})}\to 0$ and
$\|\widehat{\rho}-\widehat{\rho}_j\|_{C^{k+1,\beta}(\widehat{{U}})}
  \to 0$ as $j\to \infty$ and
\begin{eqnarray*}
(1,0,\ldots,0) & = & \left ( \frac{\partial \rho }{\partial z_1}(0),\ldots,
\frac{\partial \rho}{\partial z_n}(0) \right ) \\
& = & \left (\frac{\partial \rho_j }{\partial z_1}(0),\ldots,\frac{\partial
\rho_j }{\partial z_n}(0) \right ) \\
& = & \left (\frac{\partial \widehat{\rho }}{\partial
z_1}(0),\ldots,\frac{\partial \widehat{\rho} }{\partial z_n}(0) \right ) \\
& = &
\left (\frac{\partial \widehat{\rho }_j}{\partial
z_1}(0),\ldots,\frac{\partial \widehat{\rho}_j }{\partial
z_n}(0) \right ).
\end{eqnarray*}
\item There exist biholomorphisms $f_j \colon\Omega_j \to
\widehat{\Omega_j}, f\colon\Omega \to \widehat{\Omega}$
and a sequence $x_j \in \Omega_j$ converging to $0\in \partial
\Omega$ as $j \to \infty$
such that $f_j$ converges to $f$ uniformly on every compact
subset $K$ of $\Omega$
while $|D^\ell f_j(x_j)|\to \infty$ as $j\to \infty$ for some multi-index $\ell$
with $1\leq |\ell| \leq k$.
\end{itemize}
\vspace*{.15in}

\noindent
\bf Step 2. \it Simultaneous convexification. \rm
\medskip  \\

The content of this step is from \cite{Fornaess}.
To the expansion of  $\rho$ at $0$,
$$
\rho(z)
=2 \re  z_1+\re \sum \frac{\partial^2 \rho}{\partial
z_iz_j}(0)z_i z_j
+\frac{1}{2}\sum_{i,j}\frac{\partial^2 \rho}{\partial
z_i\overline{z}_j}(0)z_i \overline{z}_j
+o(|z|^2) \, ,
$$
apply the local biholomorphic change $\Upsilon=(w_1,w_2,...,w_n)$ of
holomorphic coordinate system at the origin $0$ defined by
$$
w_i(z)=
\begin{cases}
2 z_1+\sum \frac{\partial^2 \rho}{\partial z_iz_j}(0)z_i z_j,
& i=1 \, ,
\\
z_i,  & i =2,\ldots, n.
\end{cases}
$$
The new defining function (we continue to use $\rho$, as there is little danger of confusion) takes the form
$$
\rho
= \Re w_1
+\frac{1}{2}\sum_{i,j}\frac{\partial^2 \rho}{\partial
w_i\overline{w}_j}(0)w_i \overline{w}_j
+\varepsilon(w) \, ,
$$
where $\varepsilon(w)=o(|w|^2)$. Note that $\Upsilon(\Omega)$ is strictly convex in a small neighborhood of $0$.
Furthermore, there exists a positive integer $N$ such that  $\Upsilon(U'\cap\Omega_j)$ is strictly convex for any $j>N$.  Let
$\rho_j$ denote $\tilde{\rho}_j\circ \Upsilon$, where $\tilde{\rho}_j$is
strictly convex on $V'\cap\Omega_j$ for all $j>N$.
Set
$\displaystyle{
\tilde{\rho}(z)= \Re z_1
+\frac{1}{2}\sum_{i,j}\frac{\partial^2 \rho}{\partial
z_i\overline{z}_j}(0)z_i \overline{z}_j
+\sigma(z)
}$.
There exists a positive constant $R$ sufficiently large so that the real Hessian forms of $\tilde{\rho}(z)-\frac{|z|^2}{2R}- \Re z_1$ and $
\tilde{\rho}_j(z)-\frac{|z|^2}{2R}- \Re z_1
$
are positive-definite at every $z\in V'$.
Choose $h\in C^\infty(\RR)$ such that
$$
\begin{array}{lcr}
h(x)=0 & \hbox{if} & x \geq 1 \, , \\
0\leq h(x)\leq 1 & \hbox{if} &  0 \leq x \leq 1 \, , \\
h(x)=1 & \hbox{if} & x \leq 0 \, . \\
\end{array}
$$
Taking a larger value for $N$ if necessary, we may have that the real Hessian forms of
$$
\Re z_1+\frac{|z|^2}{2R}+\frac{1}{N}h\bigl(\frac{|z|-\eta}{\eta}\bigr)(\tilde{\rho}(z)-\frac{|z|^2}{2R}- \Re z_1)
$$
and
$$
\Re z_1+\frac{|z|^2}{2R}+\frac{1}{N}h\bigl(\frac{|z|-\eta}{\eta}\bigr)(\tilde{\rho}_j(z)-\frac{|z|^2}{2R}-\Re z_1)
$$
are both positive-definite real Hessian at every point of
$V_\delta:=\{z\in \CC^n
\colon |z|<\delta\}\subset\subset V'$ whenever
$\eta$ satisfies $0<\eta<\frac{\delta}{3}$.
Take $\eta>0$ such that $2^{2N+2}\eta<\frac{\delta}{3}$ and
set
$$
\tau(z)= \Re z_1+\frac{|z|^2}{2R}+\frac{1}{N}\sum_{m=1}^N
h\bigl(\frac{|z|-2^{2m}\eta}{2^{2m}\eta}\bigr)(\tilde{\rho}(z)-\frac{|z|^2}{2R}- \Re z_1)
$$
and
$$
\tau_j(z)= \Re z_1+\frac{|z|^2}{2R}+\frac{1}{N}\sum_{m=1}^N
h\bigl(\frac{|z|-2^{2m}\eta}{2^{2m}\eta}\bigr)(\tilde{\rho}_j(z)-\frac{|z|^2}{2R}- \Re z_1).
$$
We further let
$C=\{z\in \CC^n \colon \tau(z)<0\}$,  $C_j=\{z\in \CC^n
\colon \tau_j(z)<0\}$ and
$U^{\prime\prime}=W^{-1}(V_{\frac{\delta}{3}})$.
Then $C,C_j$ are bounded strictly convex domains
such that the restricted mappings
$
\Upsilon \big|_{U^{\prime\prime}\cap\Omega}:U^{\prime\prime}\cap\Omega \to V_{\frac{\delta}{3}}\cap
C
$ and
$
\Upsilon \big|_{U^{\prime\prime}\cap\Omega_j}:U^{\prime\prime}\cap\Omega_j \to
V_{\frac{\delta}{3}}\cap C_j
$
are biholomorphisms, and $\tau_j$ converges to $\tau$ in the $C^{k+1,\beta}$
norm, for every $\beta, 0<\beta<\alpha$.

Apply the same process to $\widehat{\Omega}$ and to
$\widehat{\Omega_j}$ at $0$.
Denote by $\widehat{C}, \widehat{C}_j$ the respective strictly convex domains
with defining functions $\widehat{\tau},\widehat{ \tau_j}$ and
$\widehat{W}:\widehat{U}\to \widehat{V}$ produced by the same procedures.
\bigskip   \\

\noindent
\bf Step 3. \it Estimates. \rm
\medskip  \\

Let $\displaystyle{\omega\in C\cap V'\cap \big(\bigcap_{j=1}^\infty C_j\big)}$ be a point
that admits an extremal map
$e : {\rm cl}(\Delta)\to {\rm cl}(C)$ satisfying
$$
e(0)=\omega, ~e(1)=0, \hbox{ and } e({\rm cl}(\Delta))\subset {\rm cl}(C)\cap V'.
$$
Let $e'(0)=\mu v$ where $|v|=1$.
Let $L : {\rm cl}(\BB^n)\to {\rm cl}({C})$ (
$L_j : {\rm cl}(\BB^n)\to {\rm cl}({C_j})$, respectively)
be the Lempert's representative map of $C$ ($C_j$, respectively) at
$\omega$.
By Proposition \ref{L}, there exists a  $\epsilon >0$ such that
$\lim_{j \to \infty} \|L_j-L\|_{C^{k,\beta}({{\rm cl}(\BB^n_\varepsilon)})}
= 0$ for any $\beta$ with $0<\beta<\alpha$.
Let $\Gamma$ be a closed cone containing $v$ in
${\rm cl}(\BB^n)$
so that $L(\Gamma)\subset {\rm cl}(C)\cap V_{\frac{\delta}{3}}$
and
$L_j(\Gamma)\subset {\rm cl}(C_j)\cap V_{\frac{\delta}{3}}$ for
all $j>N$.
Let
$$
\Upsilon^{-1}(\omega)=\zeta,
\;\;f(\zeta)=\widehat{\zeta},
\;\; f_j(\zeta)=\widehat{\zeta}_j ,
\;\;\widehat{\Upsilon}(\widehat{\zeta})=\widehat{\omega} ,
\;\;\widehat{\Upsilon}(\widehat{\zeta}_j)=\widehat{\omega}_j
$$
and let
$\widehat{L}: {\rm cl}(\BB^n) \to {\rm cl}(\widehat{C})$ and
$\widehat{L}_j: {\rm cl}(\BB^n) \to {\rm cl}(\widehat{C}_j )$, respectively, denote the Lempert representative map of $\widehat{C}$ at the point $\widehat{\omega}$ and the Lempert representative map of $\widehat{C}_j$
at the point $\widehat{\omega}_j$.

Consider now the composite maps
$\widehat{L}^{-1}\circ\widehat{\Upsilon}\circ f\circ \Upsilon^{-1}\circ L \colon
\Gamma\to \BB^n$ and
$\widehat{L_j}^{-1}\circ\widehat{\Upsilon}\circ f_j\circ \Upsilon^{-1}\circ L_j\colon \Gamma\to \BB^n$.
Denote by  $h:{\rm cl}(D)\to {\rm cl}(C) $ the extremal map satisfying
$h(0)=\omega$, $h'(0)=\lambda \zeta$, for some $\lambda >0$,
and by
$\widehat{h}=\widehat{\Upsilon}\circ f \circ \Upsilon^{-1}\circ
h:{\rm cl}(D)\to {\rm cl}(\widehat{C})$ the extremal map
satisfying
$$
\displaystyle{\widehat{h}(0)=\widehat{\omega} \ , \ \
\widehat{h}'(0)=\widehat{\lambda} |\zeta| \frac{
d(\widehat{\Upsilon}\circ f\circ \Upsilon^{-1})|_\omega (\zeta)}
{\big|d(\widehat{\Upsilon}\circ f\circ \Upsilon^{-1}|_\omega (\zeta)\big|} }
$$
for some $\widehat{\lambda}$, respectively.
Since $\widehat{C}$ is strictly convex and
$f$ extends to ${\rm cl}(\Omega)$ as $C^{k,\gamma}$
diffeomorphism for all $\gamma <\alpha$, we have
$$
\widehat{W}^{-1}\circ\widehat{L} \Bigl( \frac{|\zeta|
d(\widehat{W}\circ f\circ W^{-1})|_\omega (\zeta)}
{\big|d(\widehat{W}\circ f\circ W^{-1}|_\omega
(\zeta)\big|}\Bigr)=f\circ W^{-1}\circ L(\zeta) \, .
$$

By the same reasoning we also have
$$
\widehat{W}^{-1}_j\circ\widehat{L} \Bigl( \frac{|\zeta|
d(\widehat{W}\circ f_j\circ W^{-1})|_\omega (\zeta)}
{\big|d(\widehat{W}\circ f_j\circ W^{-1}|_\omega
(\zeta)\big|}\Bigr)=f_j\circ W^{-1}\circ L_j(\zeta).
$$
Considering the left-hand sides of the preceding identities, for any $\beta, 0<\beta<\alpha$ we obtain
$$
\lim_{j\to \infty}\|f\circ W^{-1}\circ L-f_j\circ W^{-1}\circ
L_j\|_{C^{k,\beta}(\Gamma_\varepsilon)} = 0 \, ,
$$
where $\Gamma_\varepsilon=\Gamma \setminus \lbrace z\in \Gamma \colon |z|<\varepsilon \rbrace$.
Therefore
\begin{multline*}
\lim_{j\to \infty}\|f\circ \Upsilon^{-1}\circ L-f\circ F_j^{-1}\circ
\Upsilon^{-1}\circ L_j\|_{C^{k,\beta}(\Gamma_\varepsilon)}
\\
=\lim_{j\to \infty}\|f\circ F_j^{-1}\circ \Upsilon^{-1}\circ
L_j-f_j\circ \Upsilon^{-1}\circ L_j\|_{C^{k,\beta}(\Gamma_\varepsilon)}.
\end{multline*}
Hence
\begin{eqnarray*}
\|f\circ  \Upsilon^{-1}\circ L & - & f\circ F_j^{-1}\circ \Upsilon^{-1}\circ
L_j\|_{C^{k,\beta}(\Gamma_\varepsilon)}
\\
& \lessapprox &\|\Upsilon^{-1}\circ L-F_j^{-1}\circ \Upsilon^{-1}\circ
L_j\|_{C^{k,\beta}(\Gamma_\varepsilon)}
\\
& \lesssim &\|\Upsilon^{-1}\circ L-\Upsilon^{-1}\circ
L_j\|_{C^{k,\beta}(\Gamma_\varepsilon)}
\\
&   &  \qquad +\|\Upsilon^{-1}\circ L_j-F_j^{-1}\circ \Upsilon^{-1}\circ
L_j\|_{C^{k,\beta}(\Gamma_\varepsilon)}
\\
& \lesssim    & \| L-L_j\|_{C^{k,\beta}(\Gamma_\varepsilon)}
+\|(\hbox{\rm id}-F_j^{-1})\Upsilon^{-1}\circ L_j\|_{C^{k,\beta}(\Gamma_\varepsilon)}
\\
& \to & 0 \hbox{ as } j\to\infty.
\end{eqnarray*}
On the other hand, by the proof-argument of Lemma \ref{T}, it holds that
\begin{eqnarray*}
\|f\circ F_j^{-1}\circ \Upsilon^{-1}\circ L_j & - &f_j\circ \Upsilon^{-1}\circ
L_j\|_{C^{k,\beta} (\Gamma_\varepsilon)}
\\
&=& \|(f-f_j\circ F_j)\circ F_j^{-1}\circ \Upsilon^{-1}\circ
L_j\|_{C^{k,\beta}(\Gamma_\varepsilon)}
\\
& \gtrapprox & \|(f-f_j\circ F_j)\|_{C^{k,\beta}(\sigma)}.
\end{eqnarray*}
on a sufficiently small neighborhood $\sigma$ of $p$.  This contradicts (1).
\bigskip   \\

To complete the proof let us now suppose that (2) holds.  If  $|x_j-y_j| >\kappa$ for some positive constant $\kappa$ , then
$$
\frac{|D^\gamma f_j(x_j) -D^\gamma f_j(y_j)|}{|x_j-y_j|^\beta} <
\frac{2C}{\kappa^\beta}
$$
holds for some constant $C$.  Without loss of generality,
we may assume that $x_j \rightarrow p \in \partial \Omega$
and $|x_j-y_j| < \kappa$.   Suppose that there exist sequences $x_j, y_j \in
\Omega_j$ and a positive constant $\nu$ such that $x_j\to 0 \in
\partial\Omega$ as $j\to \infty$
and $|x_j-y_j|<\nu$ so that
$$
\frac{|D^\ell f_j(x_j)-D^\ell f_j(y_j)|}{|x_j-y_j|^\beta}\to \infty
$$
as $j\to\infty$ for some multi-index $\ell $ where $|\ell |=k$.  Repeating Steps 1,2 and 3 above, we again arrive at a contradiction.  Hence the proof of
Lemma \ref{TT} is complete.
\endpf
\medskip \\

Now we present
\medskip  \\

\noindent\bf Proof of Theorem \ref{SET}. \rm
Throughout the proof, we shall take subsequences from the $\{f_j\}$ several times. But we denote them by the same notation $f_j$, since there is little danger of any confusion.

By Cauchy estimates and the standard normal family theorem, for any compact subset $K$ of  $\Omega$ we have
 $$
\lim_{j\to \infty} \|f_j-f\|_{C^{k,\beta}(K)} = 0.
 $$
Denote by  $K_\eta=\lbrace z\in \Omega \mid \hbox{dist }(\partial\Omega,z)\geq
\eta\rbrace $.
Then there exist $N>0$ and $\eta >0$ such that
$F_j(K)\subset\subset K_\eta\subset\subset \Omega$ for all
$j>N$.
So
$$
\|f_j\circ F_j-f\|_{C^{k,\beta}(K)}
\leq \|f_j\circ
F_j-f_j\|_{C^{k,\beta}(K)}+\|f_j-f\|_{C^{k,\beta}(K)} \to 0
$$
as $j\to \infty$ for all $\beta<\alpha$ by the proof of Lemma
\ref{tt}.

Let $\lambda >0$.
For $x\in {\rm cl}(\Omega)-K_\epsilon$, there exists $y\in
K_\epsilon$ such that $|x-y|<\epsilon$.
By Lemma \ref{TT}, we have
\begin{eqnarray*}
|D^l(f_j\circ F_j)(x)-D^l f(x)|
 & \leq & |D^l(f_j\circ F_j)(x)-D^l(f_j\circ F_j)(y)|
\\
&  & \ +|D^l(f_j\circ F_j)(y)-D^lf(y)| +|D^lf(y)-D^lf(x)|
\\
&  \lesssim & 2|x-y|^\beta +\epsilon
\\
& \lesssim & 2\epsilon^\beta + \epsilon .
\end{eqnarray*}
Since
\begin{multline*}
\sup_{\small x\in  {\rm cl}(\Omega) \atop 0\leq |\ell| \leq k}
|D^\ell(f_j\circ F_j)(x)-D^\ell f(x)| \\
\leq \max \Big\{\sup_{\small x\in K_\epsilon \atop \small 0\leq |\ell|
\leq k}|D^\ell(f_j\circ F_j)(x)-D^\ell f(x)|,
\sup_{x\in {\rm cl}(\Omega) \setminus K_\epsilon \atop 0\leq |\ell|
\leq k}|D^\ell(f_j\circ F_j)(x)-D^\ell f(x)|\Big\},
\end{multline*}
there exist $N>0$ and $\epsilon$ such that
$$
\sup_{x\in {\rm cl}(\Omega) \atop 0\leq |\ell| \leq k}
|D^\ell(f_j\circ F_j)(x)-D^\ell f(x)|<\lambda
$$
for all $j>N$.

Let $\displaystyle{ \delta_\ell(x,y) :=\frac{|D^\ell(f_j\circ
F_j)(x)-D^\ell f(x)-D^\ell (f_j\circ F_j)(y)+D^\ell f(y)|}{|x-y|^\beta} }$.
Then
\begin{equation}
\sup_{x,y\in{\rm cl}(\Omega) \atop |\ell|=k}
\delta_\ell (x,y)
\leq \max \biggl (\sup_{x\in {\rm cl}(\Omega) \atop y\in K_\epsilon, |\ell|=k} \delta_\ell (x,y), \sup_{x,y\in{\rm cl}(\Omega) \setminus K_\epsilon \atop |\ell|=k}
\delta_\ell(x,y) \biggr ) \, . \label{formula}
\end{equation}

Consider the first supremum in the right-hand side of (\ref{formula}).
For $x\in {\rm cl}(\Omega),y\in K_\epsilon$,
there exists $z\in K_\epsilon$ such that
$\hbox{dist }(K_\epsilon,x)=|x-z|$.
Therefore we see that
\begin{eqnarray*}
\delta_\ell (x,y)
& \leq & \frac{|D^\ell (f_j\circ F_j)(x)-D^\ell f(x) - D^\ell (f_j\circ
F_j)(z)+D^\ell f(z)|}{|x-y|^\beta}
\\
&  & +\frac{|D^\ell (f_j\circ F_j)(z)-D^\ell f(z)-D^\ell (f_j\circ
F_j)(y)+D^\ell f(y)|}{|x-y|^\beta}
\\
& \lesssim & \delta_\ell (x,z)+\delta_\ell (z,y),
\end{eqnarray*}
because $|x-y| \geq |x-z|$ and $|y-z| \leq |y-x|+|x-z| \leq
2|x-y|$.
Notice now that, for $\mu$ satisfying $\beta+\mu<\alpha$, we have that
$\delta_\ell (x,z)\lesssim |x-z|^\mu<\epsilon^\mu$.  So
$$
\sup_{x\in {\rm cl}(\Omega),y\in K_\epsilon \atop
|\ell|=k} \delta_\ell(x,y)<\lambda
$$
for any $j>N$. (For this last, one needs to adjust the sizes of $N$ and $\epsilon$
if necessary.)

Consider now the second supremum in the right-hand side of
(\ref{formula}).  Let $x, y\in{\rm cl}(\Omega)-K_\epsilon$.
If $|x-y|<\epsilon$, then for $\mu$ satisfying
$\beta+\mu<\alpha, \delta_\ell(x,y)\lesssim
|x-y|^\mu<\epsilon^\mu$.
If $|x-y| \geq \epsilon$, let $z$ be a point in $K_\epsilon$
satisfying $|x-z|= \hbox{dist }(K_\epsilon,x)$.
Then $\delta_\ell(x,y) \lesssim \delta_\ell (x,z)+\delta_\ell(z,y)$,
since $|x-z|<\epsilon<|x-y|$ and $|y-z|<2|x-y|$.
So
$$
\sup_{x,y\in{\rm cl}(\Omega) \setminus K_\epsilon \atop |\ell|=k}
\delta_l(x,y)<\lambda \, .
$$
Since $\lambda > 0$ is arbitrary, we see that
$$
\lim_{j\to\infty} \sup_{\small x,y\in{\rm cl}(\Omega) \atop |\ell|=k}
\delta_\ell (x,y) = 0
$$
for any $\beta<\alpha$.  This completes the proof of Theorem \ref{SET}.
\endpf
\smallskip \\


\section{Conjugation by Diffeomorphism}

For isometries of compact Riemannian manifolds, semicontinuity involves not just
that nearby metrics have isometry groups which are isomorphic to subgroups of the
unperturbed metric, but that the isomorphisms are obtainable via conjugation by
diffeomorphism (cf.\cite{Ebin} and \cite{Guillemin}).  This conjugation by
diffeomorphism actually applies in the case of bounded $C^\infty$ strongly
 pseudoconvex domains as well, e.g.\ \cite{Greene/Krantz_Advances, GKK}.
Naturally, the $C^\infty$ hypothesis used in these references is, as usually happens,
replaceable by a finite differentiability hypotheses simply by tracing through the
arguments and checking how many derivatives are needed.
\medskip  \\

In this section, the subject will be investigated of  the finite differentiability version of
the conjugation by diffeomorphism results already shown in the references indicated in
the $C^\infty$ case.  These results are of active interest because, by this time, quite
precise results are known about extension to the boundary with finite smoothness of
automorphisms of bounded  strongly pseudoconvex domains with boundaries of finite
smoothness. In particular, the results of the previous sections give motivation to study
the issues discussed in the present section.
\medskip  \\

In the $C^\infty$ version presented in \cite{Greene/Krantz1982} and \cite{GKK}, the
basic technique was to pass to the double in the topologist's sense of the domain, thus
creating a situation to which the compact manifold results could be applied.  This
technique can still be applied in the present case. The difference is that we need now to
keep track of how many derivatives are lost in the passage to the double. For the
manifold with boundary itself, no derivatives are lost. It is shown in \cite{Munkres}
that a $C^k$ manifold with boundary, $k \ge 1$, has a $C^k$ double that is unique
up to $C^k$ diffeomorphism.

But, in our case, the need to make the group act on the double requires that the
doubling construction be invariant under the group. And this will turn out to
reduce the guaranteed differentiability of the conjugating diffeomorphism.

To facilitate the discussion, we introduce a definition (similar to one given in Section 2)
of the sense in which a
sequence of groups of diffeomorphisms might converge to a limit group:

Suppose $M$ is a compact $C^k$ manifold with boundary, $k$ a positive integer.
Suppose $G_0$ is a compact Lie group of $C^k$ diffeomorphisms of $M$ and that
moreover  $G_j$,  $j=1,2, \ldots$ are a sequence of compact Lie groups of $C^k$
diffeomorphisms.
Then we say that the sequence $G_j$ converges to $G_0$ in the $C^k$ sense if for each
$\epsilon >0$ there is a number $j_0$ such that, if $j >j_0$ and $g \in G_j$ , then there is
an element $g_0 \in G_0$ such that the distance from $g$ to $g_0$ is less than
$\epsilon$. Here the distance means relative to any metric on the set of $C^k$
mappings which gives the usual $C^k$ topology on $C^k$ maps from $M$ to $M$.

In these terms, we can now formulate the general real-differentiable result we shall use
in the complex case:

\begin{theorem} \label{conjugation-th}  \sl
Suppose that $M$ is a compact $C^3$ manifold with boundary and that $r>2$ is an integer, that $G_0$ is a compact Lie group of $C^r$ diffeomorphisms of $M$ , and that $G_j$, $j=1,2,\ldots$, is a sequence of compact groups of $C^3$ diffeomorphisms which converge in the $C^3$ sense to $G_0$. Then, for all $j$ sufficiently large there is a $C^{r-2}$ diffeomorphism $F_j$ of $M$ to itself such that $F_j \circ G_j\circ F_j^{-1}$ is a subgroup of $G_0$, i.e., $F_j$ conjugates the elements of $G_j$ into elements of $G_0$.
\end{theorem}

The proof of this theorem follows almost precisely the pattern of the proof of
Theorem 0.1 in \cite{Greene/Krantz1982} (cf.\ Theorem 4.4.1 \cite{GKK}). The only
difference is that we must here keep some track of the number of derivatives involved:
Ebin's theorem concerned the $C^\infty$ case so that loss of a derivative or two or
indeed of any finite number was irrelevant.
\medskip  \\

As in Section 2, the essential method is to pass to the double of $M$ and
extend the action of the groups
to the double. Then one can use Ebin's result in the form presented in \cite{Guillemin},
where only $C^1$ is required for the closeness of the group actions.
But here we have to keep track of degrees of differentiability as opposed to
the $C^\infty$ situation of Section 2.

The most natural way to form the equivariant double is via metric construction as
already explained in Section 2 (cf.\ \cite{Greene/Krantz1982}): As before one
takes a metric $g$ on the manifold with boundary that is
invariant under the group $G$. Then one defines charts in neighborhoods of boundary
points $p$ using the normal field to the boundary. Specifically, let
$N(q)$ be the $g$-metric normal to the tangent space to the boundary $\partial M$ at
the point $q$ in $\partial M$.  Then  one defines charts in a neighborhood of points $p$
in the boundary as follows: Map $\partial M \times (-\epsilon, \epsilon) \to M$  by
$(q, t) \mapsto \exp_q (t N(q))$, where $\exp$ is the geodesic exponential map of the
Riemannian metric $g$ and $N(q)$ is the inward pointing normal at $q$.   Choosing a
chart around $p$ in $\partial M$ then gives a chart in a neighborhood of $p$ in the
double of $M$ if we interpret $\exp_q (tN(q))$ to be in the second copy of $M$
when $t <0$.

In terms of derivative loss, the choice of the normal vector $N$ loses one derivative, since it is an
algebraic process using $g$ and the tangent space of the boundary and the latter is not
$C^r$ but $C^{r-1}$.   But an additional loss of derivative, so that two derivatives are
lost, occurs because the exponential map is defined by the geodesic equation and that
equation involves the Christoffel  symbols, which involves the first derivative of the
metric $g$. And the metric $g$ has already lost one derivative in the averaging over the
action of the group $G$.

Thus one obtains a $G$-equivariant construction of the double $\widetilde M$ of $M$
and by construction the action of $G$ on $M$ extends to be an action of $G$ on
$\widetilde M$. This extended group action is $C^{r-2}$.   Associate to the group $G$ a
group $\widetilde G$ defined to be $G \oplus \ZZ_2$. Then $\widetilde G$
acts on $\widetilde M$ in a natural way. Namely, we label the elements of $\widetilde
M$ by $(m, a)$ where $m \in M$ and $a \in  \{0, 1\}$ with $0$ corresponding to the
original of $M$ and $1$ corresponding to the second copy of $M$. Then we let $(g,b)$
acting on $(m,a)$ be $(g(m), a+b)$ where the addition $a+b$ is in $\ZZ_2$, i.e.,
mod $2$.  For example $(\hbox{id}_G, 1)$ acts on $\widetilde M$ as the ``flip'' map that
interchanges the two copies of $M$.

Note that the fixed point set of $(\hbox{id}_G,1)$ is exactly $\partial M$. And, for any
element $g \in G$, the fixed point set of $(g,1)$ is contained in $\partial M$, though it
need not be all of it, and can indeed be empty if the action of $g$ on $\partial M$ has no
fixed point. These observations will be important later.

Now we turn to the explicit situation of Theorem \ref{conjugation-th}. We choose a
sequence of $G_j$-invariant $C^{r-1}$  metrics on $M$, which can clearly be taken to
converge in the $C^{r-1}$ sense to a $G_0$-invariant $C^{r-1}$ metric on $M$.   Passing to
the double $\widetilde M$ gives a sequence of $\widetilde G_j$ group actions on
$\widetilde M$. We can form a sequence of $\widetilde G_j$ invariant metrics by
combining, via a partition of unity, a product metric structure near the boundary with
the $G_j$-invariant metric on the interior of $M$. Namely, as similar to before, let $E_j$
be the exponential map of the metric $g_j$, $j=0,1,2,\ldots$, acting on the normal bundle
of the boundary $\partial M$ of $M$ in $M$ to give maps also to be denoted by $E_j
\colon \partial M \times (-a,a) \to \widetilde M$
of the boundary $\partial M$ of $M$ producted with an open interval $(-a,a)$ into
$\widetilde M$. The size of $a$ can, by the $C^{r-2}$ convergence of the $E_j$ to $E_0$,
be chosen uniformly so that these $E_j$ are diffeomorphisms onto their images in
$\partial M$, which themselves converge in the $C^{r-2}$ sense to the limit  $C^{r-2}$
diffeomorphism $E_0$.

Via this diffeomorphism, we transfer the product metrics on $\partial M \times [0,
\epsilon)$, namely $H_j \times dt^2$, to the associated tubular neighborhoods of
$\partial M$ in $M$.  This transfer gives a $\widetilde G_j$-invariant metric for each $j$
and these metrics converge $C^{r-2}$ to the limit $\widetilde G_0$-invariant metric.
Now we can combine, using a $\widetilde G_j$-equivariant partition of unity, these
product metrics with the $G_j$-invariant metric $g_j$ on $M$ to obtain a $\widetilde
G_j$-invariant metric on $\widetilde M$, to be denoted $\widetilde g_j$. This metric is
$C^{r-2}$. And it converges in the $C^{r-2}$ topology to the corresponding $\widetilde
G_0$-invariant metric $\widetilde g_0$ on $\widetilde M$.   (The $G_j$-equivariant
partition of unity is obtained by taking the partition of unity function to depend on $t$
alone, $t$ as above).

Now we can apply Ebin's Theorem, in the form given in \cite{Guillemin} and
\cite{KimY}, for the $C^{r-2}$ case to get $C^{r-2}$ diffeomorphisms $F_j\colon
\widetilde M \to \widetilde M$ which conjugate $\widetilde  G_j$ into a subgroup of
$\widetilde  G_0$.  (Here we are reasoning as follows: There is a diffeomorphism that
conjugates $\hbox{Isom }(\widetilde g_j)$ into a subgroup of $\hbox{Isom }
(\widetilde g_0)$ and hence conjugates $\widetilde G_j$ into a subgroup of
$\hbox{Isom } (\widetilde g_0)$ and these diffeomorphisms can be taken to converge to
the identity map. So the image of $\widetilde  G_j$ under this conjugation is close to
$\widetilde  G_0$ for large $j$ in the sense of $C^{r-2}$ convergence. By the classical
theorem of \cite{Montgomery/Samelson}, this conjugation image is in fact itself
conjugate in $\hbox{Isom }(\widetilde g_0)$ to a subgroup of $G_0$ by an element
close to the identity. (cf., e.g., \cite{GKK}, Ch.\ 4, for more detail.)

Now we need to know that in fact the conjugation image of $G_j$ lies in $G_0$,
not just in $\widetilde G_0$. For this, we need only show that the diffeomorphism
that is conjugating takes $\partial M$ to itself.  This can be deduced as follows: Let us
denote by $\hbox{Fix }(\psi)$ the fixed point set of $\psi$. Then conjugation takes fixed
points to fixed points in the sense that
$\hbox{Fix }(f \circ \psi \circ f^{-1}) = f(\hbox{Fix }(\psi))$.
Now consider the case of $\psi=$ the flip map which interchanges the two copies of
$M$ in $\widetilde M$. When $f$ is close to the identity,  $f\circ \psi\circ f^{-1}$ has to
belong to the part of the group that interchanges the two components. So its fixed point
set cannot be larger than $\partial M$.  Thus $f(\partial M)$ lies in $\partial M$ and
hence equals $\partial M$ (since $f$ is a diffeomorphism of $\partial M$ onto its
image).

This completes the proof of the theorem. \endpf
\bigskip  \\

Note that these considerations of fixed points of the interchange map did not arise in Section 2, since we were concerned there only with isomorphism, not with the existence of a conjugating diffeomorphism of the manifolds with boundary.

The application to the strongly pseudoconvex case now follows:

\begin{theorem}  \sl
Let $\Omega_0$ be a bounded strongly pseudoconvex domain with a $C^{4,\alpha}$ boundary in $\CC^n$, not biholomorphic to the unit ball. Then there is a $C^{4,\alpha}$ neighborhood $\mathcal N$ of $\Omega_0$ such that, for any $\Omega \in \mathcal N$, there is a $C^3$ diffeomorphism $f \colon \Omega \to \Omega_0$ with the property that $f \circ \Aut (\Omega) \circ f^{-1} \subset \Aut (\Omega_0)$.
\end{theorem}
\vfill
\eject

\null 
\vspace*{.3in}

\begin{quote}
(Greene) Department of Mathematics, University of California , Los Angeles, California 90095-1555, U. S. A. \ \ \ {\tt greene@math.ucla.edu} \\
\smallskip   \\
(Kim) Department of Mathematics, Pohang University of Science and Technology, Pohang 790-784 South Korea \ \ \ {\tt kimkt@postech.edu}  \\
\smallskip \\
(Krantz) Department of Mathematics, Washington University in St. Louis, St. Louis, Missouri 63130-4899 U. S. A. \ \ \ {\tt sk@math.wustl.edu}  \\
\smallskip \\
(Seo) Department of Mathematics, Pohang University of Science and Technology, Pohang 790-784 South Korea. \ \ \ {\tt inno827@postech.edu}
\end{quote}

\end{document}